\documentclass[a4paper]{article}

\usepackage{latexsym}

\usepackage{amsmath} 
\usepackage{amssymb} 
\usepackage{amsthm}  

\newtheorem{thm}{Theorem}[section]

\newtheorem{lem}[thm]{Lemma}
\newtheorem{cor}[thm]{Corollary}

\def\qed{{$\Box$}\medskip}
\def\altem{\slshape}

\def\a{\alpha} \def\b{\beta} \def\g{\gamma} \def\d{\delta}
\def\e{\varepsilon}  \def\k{\kappa} \def\l{\lambda}
 \def\r{\rho} \def\s{\sigma} \def\ta{\theta} \def\t{\tau}
\def\om{\omega} \def\Om{\Omega} \def\dd{\partial} \def\D{\Delta}
\def\G{\Gamma}

\def\bC{{\mathbf C}} \def\bN{{\mathbf N}} 
\def\bR{{\mathbf R}}   \def\bZ{{\mathbf Z}}

 \def\id{\mathrm{id}}

\def\ip#1#2{{(\,#1\mid #2\,)}}
\def\m#1{{\vert #1 \vert}} \def\n#1{{\Vert #1 \Vert}}

\def\demo#1{\noindent \emph{#1.}}

\def\sa{\mathcal{A}} \def\sag{\sa_G} \def\ominv{\Om^{\mathrm{inv}}}
\def\omhar{\Om_\nabla}
\def\omd{{(\Om,d)}}
\def\E{{\mathcal E}} \def\GG{{\mathcal G}}
\def\EEONE{{\E_1(M,p,k)}}
\def\EE{{\E(M,p,k)}}
\def\F{{{\mathcal F}_1}} \def\FF{{\F(M,p,k)}}
\def\T{{\mathcal T}}

\def\N{\bN}
\def\Z{\bZ}
\def\R{\bR}
\def\C{\bC}

\def\spat{\hspace{4ex}}
\def\flip{\psi}

\def\nab{\nabla}

\def\ot{\otimes}

\begin{document}

\title{Quantum groups, differential calculi and the eigenvalues of
the Laplacian}
\author{J. Kustermans \and G.J. Murphy \and L. Tuset}
\date{}

\maketitle

\begin{abstract}  We study $*$-differential calculi
over compact quantum groups in the sense of S.L.~Woronowicz. Our principal results are the construction of a Hodge
operator commuting with the Laplacian, the derivation of a corresponding Hodge decomposition of the calculus of forms,
and, for Woronowicz' first calculus, the calculation of the eigenvalues of the Laplacian. 
\end{abstract}

\section{Introduction}

 Compact quantum groups and differential calculi over
them were introduced by S.L.~Woronowicz in a series of seminal papers~\cite{WOR1,WOR2,WOR3}. The theory of compact
quantum groups is now well established---we refer to~~\cite{KT} and~\cite{MT} for recent expositions. On the other
hand, the theory of differential calculi over quantum groups is still in a very rudimentary state of development. A
possible reason for this is that although it is an important part of non-commutative geometry, it does not fit into the
framework of A.~Connes. The principal difference from Connes theory is that the appropriate ``integrals'' on quantum
differential calculi are not graded traces but twisted graded traces. This involves the appearance of a twist
automorphism under the integral sign and the transition from graded traces to twisted graded traces is analogous to
that from ordinary traces in the theory of von~Neumann algebras to KMS~states. Some of the implications of this
transition, such as the necessity of developing a twisted cyclic cohomology theory, are explored in~\cite{KMT}.

In this paper we are principally concerned with studying the
analogues of the classical operators of ordinary differential
calculus over a manifold. Thus, we show that for the calculi that
we study the differential admits an adjoint or co-differential and
therefore one can define Dirac and Laplace operators. We also show
that with suitable natural hypotheses on the calculi a Hodge
operator exists. We then show that the calculus of forms admits a
Hodge decomposition. This is surprising, since such a
decomposition is known not to exist in general in the setting of
non-commutative differential geometry.

We derive some general results on the eigenvalues of the
Laplacian; specifically, we obtain formulas for the
eigenvalues of the Laplacian when restricted to the spaces of
lowest and highest dimensional forms. (The eigenvalues of the
Dirac operator can be calculated from those of the Laplacian.)

The differential calculi that we are concerned with are
assumed to be left invariant (this involves no loss of
significant generality, since all the interesting examples in
the theory are left invariant).  The choice of a suitable
inner product on the space of forms reduces to the question of
choosing a suitable inner product on the space of invariant
forms, which is finite-dimensional in our theory. For
many results this choice is quite arbitrary, since one needs
only assume that the inner product is chosen so that forms of
different dimensions are orthogonal. However, we believe that
for the further development of the theory the choice of inner
product and what kind of properties it should have will play a
more significant role. We take up this question and show that
one can always choose the inner product in such a way that the
Hodge operator  corresponding to this choice has such
desirable properties as being a unitary that commutes with the
Laplacian. We also show that the inner product can be chosen
so that the twist automorphism corresponding to the ``volume
integral'' on the calculus of forms is positive.

Connected with this, we show that under  a natural condition on the twist automorphism, this twist automorphism can be
seen as the generator of a one-parameter group of automorphisms on the space of differential forms. This one-parameter
group mimics the role played by the modular groups of normal faithful semi-finite weights on von Neumann algebras. It
also  provides us with a symmetry group like the one discussed in \cite[VII.12]{Jaffe}.

The important question of calculating {\em all} the
eigenvalues of the Laplacian is undertaken in the case of a
specific calculus, namely the three-dimensional calculus of
Woronowicz over quantum $SU(2)$ that was the first ever
example of a quantum calculus. The formulas for the
eigenvalues are complicated, but we are nevertheless able to
obtain bounds for them that enable us to show a ``scaled
boundedness'' result for the commutator of the Dirac operator
and an arbitrary element of the Hopf algebra underlying the
quantum group. This shows that although this differential
calculus does not fit into Connes' framework, it has
properties that are close to those required to enable it to
fit into the framework developed by A.~Jaffe (this theory is
an extension of Connes' theory). It may be that Jaffe's theory
can, and should, be further developed to cover this example
(and the quantum differential theory in general)---in any
case, our example illustrates what requires to done.

A very nice solution to similar kinds of problems has been developed in \cite{Heck1} and \cite{Heck2} by
I.~Heckenberger and A.~Sch\"{u}ler. In these two papers, the authors look at metrics which they assume to satisfy the
locality property but as a result the considered metrics cannot be positive in some cases (like in the case of
Woronowicz' 4-dimensional differential calculi on quantum $SU(2)$\,). Nevertheless, a powerful theory emerges in their
setting.

Commutator representations of differential calculi and the unboundedness of them are also intensively studied by
K.~Schm\"udgen, see for instance \cite{KS1} and \cite{KS2}.

We summarise now what we do in each section of the paper:
Section~2 introduces the basic concepts and terminology that
will be needed throughout the paper and shows the existence of
the co-differential under suitable hypotheses
(Theorem~\ref{thm: d-adjointability theorem}). Section~3 is
concerned with the existence of the Hodge operator
(Theorem~\ref{thm: existence of Hodge operator}) and in
Section~4 a Hodge decomposition is obtained (Theorem~\ref{thm:
Hodge decomposition}). In Section~5 we prove some general
results on diagonalising the Laplacian on the forms of lowest
and highest degree. Sections~6 and~7 are concerned with
choosing the inner product on the invariant forms in such a
way that the Hodge operator and the twist operator have
desirable properties (Theorems~\ref{thm: L-properties theorem}
and~\ref{thm: inner product theorem}). Section~8 is devoted to
the calculation and analysis of the eigenvalues of the
Laplacian for the three-dimensional calculus of Woronowicz
mentioned above. At the end of the paper we have included a
brief appendix on some non-standard results on exponential
one-parameter automorphism groups that we need.

\section{Differential calculi}
\label{sec: differential calculi}

We begin by recalling some basic definitions.

If $\Om$ is an algebra, a {\em positive algebra grading} on
$\Om$ is a sequence $(\Om_n)_n$ of subspaces of $\Om$ for
which $\Om=\oplus_{n=0}^\infty \Om_n$ and $\Om_m\Om_n\subseteq
\Om_{n+m}$, for all $n,m\ge 0$. The pair $(\Om, (\Om_n)_n)$ is
called a {\em graded algebra}.

A conjugate-linear map ${\om \mapsto \om^*}$ of degree zero is
called a {\em graded involution} on $\Om$ if $(\om^*)^*=\om$
and $(\om_1\om_2)^*=(-1)^{kl}\om_2^*\om_1^*$, for all $\om\in
\Om$ and all pairs $\om_1,\om_2\in \Om$ of degrees $k$ and
$l$, respectively. The pair consisting of $\Om$ together with
this involution is called a {\em graded $*$-algebra}.

A {\em graded $*$-derivation} on $\Om$ is a linear map
${d\colon \Om\to \Om}$ for which $d(\om^*)=(d\om)^*$ and
$d(\om'\om)=d(\om')\om+(-1)^k\om' d\om$, for all $\om'\in
\Om_k$ and all $\om\in \Om$.

A {\em graded $*$-differential algebra} is a pair $(\Om,d)$,
where $\Om$ is a graded $*$-algebra, $d$ is a graded
$*$-derivation on $\Om$ of degree 1 (as a linear map) and
$d^2=0$. The elements of $\Om$ are referred to as the {\em
forms} of $(\Om,d)$ and the elements of $\Om_k$ as the {\em
$k$-forms}. The operator $d$ is referred to as the {\em
differential}.

Suppose now $(\sa,\D)$ is a Hopf $*$-algebra. (We shall
usually refer to the algebra by the symbol $\sa$ and suppress
explicit reference to the co-multiplication~$\D$. For
Hopf algebra theory we refer the reader
to~\cite{ABE}.) We say that a graded $*$-differential algebra
$(\Om,d)$ is a {\em $*$-differential calculus over $\sa$} if
the following conditions hold:

\begin{enumerate}

\item $\Om_0=\sa$ (as $*$-algebras) and, for all $k>0$, $\Om_k$
is the linear span of all products $a_0da_1\cdots da_k$, where
$a_0,\dots,a_k\in \sa$;

\item There exists an algebra homomorphism ${\D_\Om\colon \Om\to
\sa\otimes \Om}$ such that $\D_\Om(a)\\=\D(a)$, for all $a\in
\sa$, and such that ${(\id_\sa\otimes d)\D_\Om}=\D_\Om d$;

\item d1=0.

\end{enumerate}

Note that Condition~2 is called {\em left covariance} of the
differential calculus. (We shall only be interested in
left-covariant differential calculi.)

We now make a sequence of useful observations that follow from
this definition:

Conditions~1 and~3 easily imply that the unit 1 of $\sa$ is a
unit for the algebra~$\Om$.

Since $d(a b)=(d a)b+adb$, for $a,b\in \sa$, it follows from
Condition~1 that $\Om_1$ is the linear span of all products
$(da_1)a_0$, where $a_0,a_1\in \sa$. A simple inductive
argument now shows that $\Om_k$ is the linear span of all
products $(da_k)\cdots (da_1)a_0$, where $a_0,\dots,a_k\in
\sa$.

It is easily seen that the map $\D_\Om$ is automatically of
degree zero and $*$-preserving; that is,
$\D_\Om(\om^*)=\D_\Om(\om)^*$.

Note that ${(\D\otimes \id_\Om)\D_\Om}={(\id_{\sa}\otimes
\D_\Om)\D_\Om}$ and ${(\e\otimes \id_\Om)\D_\Om}=\id_\Om$,
where $\e$ is the co-unit of $\sa$; hence, the map $\D_\Om$ is
a left co-action of the Hopf algebra $\sa$ on $\Om$. For the
proof, see~\cite[Lemma~4.2]{KMT}.

Recall that a bi-module $\Gamma$ over $\sa$ is said to be {\em
left covariant} if there is a linear map
${\D_\G\colon\Gamma\to\sa\otimes\Gamma}$ such that
${(\D\otimes \id_\Gamma)\D_\G}={(\id_\sa\otimes \D_\G)\D_\G}$
and ${(\e\otimes \id_\Gamma)\D_\G}=\id_\Gamma$, where $\e$ is
the co-unit of $\sa$ (that is, $\D_\G$ is a left co-action)
and $\D_\G(a\g b)=\D(a)\D_\G(\g)\D(b)$, for all $\g\in \Gamma$
and $a,b\in \sa$. An element $\g\in \Gamma$ is said to be {\em
left invariant} if $\D_\G(\g)={1\otimes \g}$. We denote by
$\Gamma^{\mathrm{inv}}$ the linear space of left-invariant
elements of $\Gamma$.

We shall make use of the following results from the theory of
left-covariant bi-modules~\cite{KS} (we do not state these
results in their strongest forms). First, we need some
notation. If $\chi$ is a linear functional on a Hopf algebra
$\sa$ and $a\in \sa$, we write $\chi * a $ for ${(\id\otimes
\chi)\D(a)}$ and $a*\chi$ for ${(\chi\otimes \id)\D(a)}$.

\begin{thm}\label{thm: tensor isomorphism} Let $\Gamma$ be a
left-covariant bi-module over a Hopf $*$-algebra $\sa$. Then
there is a unique isomorphism of left $\sa$-modules from
${\sa\otimes \Gamma^{\mathrm{inv}}}$ onto $\Gamma$ that maps
${a\otimes \g}$ onto $a\g$, for all $a\in \sa$ and $\g\in\
\Gamma^{\mathrm{inv}}$. \end{thm}

\begin{thm}\label{thm: basis theorem} Let $\sa$ be a Hopf
$*$-algebra with co-inverse $\k$. Let $\Gamma$ be a
left-covariant bi-module over $\sa$ and suppose that
$\g_0,\dots,\g_n$ is a linear basis for
$\Gamma^{\mathrm{inv}}$. Then $\g_0,\dots,\g_n$ is a free left
$\sa$-module basis for $\Gamma$ and also a free right
$\sa$-module basis of $\Gamma$. Moreover, there exist linear
functionals $f_{jk}$ on $\sa$, for ${j,k=0,\dots,n}$, such
that $f_{jk}(a b)=\sum_{r=0}^n f_{jr}(a)f_{rk}(b)$ and
$f_{jk}(1)=\d_{jk}$ and for which we have the equations
$\g_ja={\sum_{k=0}^n(f_{jk}*a)\g_k}$ and $a\g_j= {\sum_{k=0}^n
\g_k((f_{jk}\k^{-1})*a)}$.\end{thm}

If $(\Om,d)$ is a $*$-differential calculus over a Hopf
$*$-algebra $\sa$, then it is clear from our earlier
observations that $\Om$ is a left-covariant bi-module
over~$\sa$, so that Theorem~\ref{thm: tensor isomorphism}
applies. We shall identify ${\sa\otimes \ominv}$ and $\Om$, as
left $\sa$-modules, using the isomorphism of the theorem.
Observe that each space $\Om_k$ is also a left-covariant
bi-module over~$\sa$.

It is clear that $\ominv$ is not only a linear subspace of
$\Om$, but also a subalgebra closed under the involution
${\om\mapsto \om^*}$. For $k\ge 0$, set
$\ominv_k=\Om_k\cap\ominv$. It is easily seen that
$\ominv=\oplus_{k=0}^\infty \ominv_k$, so that we have an
induced algebra grading on $\ominv$; hence, $\ominv$ is a
graded $*$-algebra. Since $\ominv$ is clearly invariant
under~$d$, $(\ominv,\dd)$ is a graded $*$-differential
algebra, where $\dd$ is the restriction of $d$ to $\ominv$.

Suppose now that $\sa$ admits a {\em Haar integral} $h$; that
is, $h$ is a unital linear functional on~$\sa$ for which
${(\id\otimes h)\D(a)}= {(h\otimes \id)\D(a)}=h(a)1$, for all
${a\in \sa}$, where $1$ is the unit of~$\sa$. We suppose
further that $h$ is positive; that is, ${h(a^*a)\ge 0}$, for
all $a\in\sa$. It is well known $h$ is unique and
automatically faithful; that is, if $h(a^*a)=0$, then $a=0$.
We endow $\sa$ with an inner product derived from~$h$, by
setting $\ip ab=h(b^*a)$, for $a,b\in \sa$.

Now suppose that $(\Om,d)$ is {\em strongly finite
dimensional}; by this we mean that $\ominv$ is finite
dimensional (in this case Theorem~\ref{thm: basis theorem}
applies to $\Om$). It follows that, for some integer $N$,
$\Om_k=0$ for $k> N$ and $\Om_N\ne 0$. Hence, $\Om$ is
finite-dimensional in the sense of~\cite{KMT}. The integer $N$
is the {\em dimension} of $\omd$.

Choose an arbitrary inner product on $\ominv$ for which the
subspaces $\ominv_k$ are pairwise orthogonal. Then endow
$\Om=\sa\otimes \ominv$ with the tensor product inner product.
Since $\Om_k={\sa\otimes \ominv_k}$, the spaces $\Om_k$ are
pairwise orthogonal. We call an inner product on~$\Om$ {\em
graded} if it has this property. We call an inner product {\em
left-invariant} if $\ip{a\om}{b\eta}=h(b^*a)\ip\om\eta$, for
all $a,b\in \sa$ and $\om,\eta\in\ominv$.

Suppose that $\dim(\ominv_N)=1$. Then it is easily seen that
there is a unit vector $\ta$ of $\ominv_N$ such that
$\ta^*=\ta$ and that $\ta$ is unique up to sign. We shall call
$\ta$ and $-\ta$ the {\em volume elements} of $\omd$.

We shall now derive a very useful formula for the
differential~$d$. First we need to introduce some linear
operators. If $\chi$ is a linear functional on $\sa$, we
denote by $E_\chi$ the linear operator on $\sa$ defined by
setting $E_\chi(a)=\chi*a$. If $\om\in \ominv_1$, we denote by
$M_\om$ the linear operator on $\ominv$ defined by setting
$M_\om(\eta)=\om\eta$.

For the next theorem, let us recall the definition of the {\em
adjoint} $\chi^*$ of $\chi$: It is given by the formula
$\chi^*(a)=\overline{\chi(\k(a)^*)}$, for all $a\in \sa$. We
shall also need to define $\bar\chi$ in $\sa^*$ by
$\bar\chi(a)=\chi(a^*)^-$.

If $T$ is an operator on an inner product space $X$, we shall
denote by $T^*$ the unique operator---if it exists---for which
$\ip{T(x)}y=\ip x{T^*(y)}$, for all $x,y\in X$. In the case
$T^*$ exists, we say $T$ is {\em adjointable} and call $T^*$
the {\em adjoint} of $T$. (Of course, if $X$ is a Hilbert
space and $T$ a bounded linear operator on $X$, $T^*$
necessarily exists. The problem of existence arises with
unbounded operators and incomplete spaces.) If $T$ is
adjointable and $T=T^*$, then the eigenvalues of $T$ are real.
Similarly, the eigenvalues of $TT^*$, $T^*T$ and $T^*T+TT^*$
are non-negative. (The usual trivial proofs apply.) .

\begin{lem} \label{thm: E-adjointability theorem}Let $\sa$
be a Hopf $*$-algebra admitting a Haar integral~$h$ and let
$\chi$ be a linear functional on $\sa$. Then $E_\chi$ is
adjointable and $E^*_\chi=E_{\chi^*}$. \end{lem}


\demo{Proof} Let $a,b\in \sa$. Using the equality $\kappa((h \otimes \id)((q \ot 1)\D(p)))=$ \newline $(h \otimes
\id)(\D(q)(p \ot 1))$ (see the proof of proposition~3.11 of \cite{VD}), we see that
\begin{eqnarray*}
&(\,E_\chi (a)\mid b\,) = h(b^*(\id\otimes\chi)\D(a))=
\chi((h\otimes \id)((b^*\otimes 1)\D(a)))\\
&=\chi(\k^{-1}((h\otimes \id)(\D(b^*)(a\otimes 1))))=
h((\id\otimes \chi\k^{-1})\D(b^*)a)\\ &= h(((\id\otimes
\chi^*)\D(b))^*a)= h(E_{\chi^*}(b)^*a)=\ip a{E_{\chi^*}(b)}.
\end{eqnarray*}
Hence, $E_\chi$ has adjoint $E_{\chi^*}$, as required. \qed

If $\om_1,\dots,\om_M$ is an orthonormal basis for $\ominv_1$,
there exist unique linear functionals
${\chi_1,\dots,\chi_M}$ on $\sa$ such that
\begin{equation} \label{eqn: d-chi equation} d a={\sum_{r=1}^M
(\chi_r*a)\om_r}.\end{equation} (See~\cite{KS}.)

\begin{thm} \label{thm: d-adjointability theorem} Let $(\Om,d)$
be a strongly finite-dimensional $*$-differential calculus
over a Hopf $*$-algebra~$\sa$ admitting a Haar integral. Then
$d$ is adjointable and $d^*$ is of degree $-1$. Indeed, if
$\om_1,\dots,\om_M$ and $\chi_1,\dots,\chi_M$ are as in
Equation~(\ref{eqn: d-chi equation}), then
\begin{equation} \label{eqn: tensor d-equation}
d=\id_\sa\otimes \dd+\sum_{j=1}^M E_{\chi_j}\otimes
M_{\om_j}\end{equation} and
\begin{equation} \label{eqn: tensor d*-equation}
d^*=\id_\sa\otimes \dd^*+\sum_{j=1}^M E_{\chi_j^*}\otimes
M_{\om_j}^*.\end{equation}
\end{thm}

\demo{Proof} Recall that we identify $\Om$ and ${\sa\otimes
\ominv}$ by identifying $a\om$ with ${a\otimes \om}$, for
$a\in \sa$ and $\om\in \ominv$. Since $d(a\om)=(d
a)\om+ad\om=\sum_{j=1}^M(\chi_j*a)\om_j\om+ad\om$, we have
$d={\id_\sa\otimes \dd}+{\sum_{j=1}^M E_{\chi_j}\otimes
M_{\om_j}}$. Adjointability of $d$ now follows from
adjointability of the operators $E_{\chi_j}$, and $M_{\om_j}$
and $\dd$. The operators $E_{\chi_j}$ are adjointable by
Lemma~{\ref{thm: E-adjointability theorem}, and the other
operators are adjointable since $\ominv$ is finite
dimensional. The formula for $d^*$ in the statement of the
theorem follows immediately. Since $d$ is of degree $+1$, it
is trivially verified that $d^*$ is of degree~$-1$.~\qed

The operator $d^*$ is called the {\em co-differential} of $d$. The sum $D=d+d^*$ is the {\em Dirac operator} and the
square $\nabla=(d+d^*)^2$ the {\em Laplacian}. Since $d^2=d^{*2}=0$, we have $\nabla=dd^*+d^*d$. We call the forms
$\om\in \Om$ such that $\nabla(\om)=0$ the {\em harmonic forms} of $(\Om,d)$ and denote the linear space of these forms
by~$\omhar$.

\begin{thm}
The co-differential, the Dirac operator and the Laplacian commute with the coaction:
$$\D_\Om d^* = (\id_\sa \ot d^*)\D_\Om \hspace{6ex}  \D_\Om\,D =(\id_\sa
\ot D)\D_\Om \hspace{6ex} \D_\Om\,\nab = (\id_\sa \ot \nab)\D_\Om \ .$$
\end{thm}
\demo{Proof} Define the linear operator $W : \sa \ot \Om \rightarrow \sa \ot \Om$ so that $W(a \ot \om) = \D_\Om(\om)(a
\ot 1)$, for all $\om \in \Om$ and $a \in \sa$. On $\sa \ot \Om$, we consider the tensor product of the natural inner
products on $\sa$ and $\Om$. For $a_1,a_2,b_1,b_2 \in \sa$ and $\om_1,\om_2  \in \ominv$,
\begin{eqnarray*}
& & \ip{W(a_1 \ot b_1 \om_1)}{W(a_2 \ot b_2 \om_2)} = \ip{\D(b_1)(a_1 \ot \om_1)}{\D(b_2)(a_2 \ot \om_2)}
\\ & & \hspace{1ex} = (h \ot h)((a_2^* \ot 1)\D(b_2^* b_1)(a_1 \ot 1))\, \ip{\om_1}{\om_2}
\\ & & \hspace{1ex} = h(a_2^* a_1) \,h(b_2^*b_1) \, \ip{\om_1}{\om_2} = \ip{a_1 \ot b_1 \om_1}{a_2 \ot b_2 \om_2} \ ,
\end{eqnarray*}
where we used the left invariance of $h$ in the second equality. This implies that $W$ is isometric. By the remarks
after \cite[Thm~4.1]{KMT} we know that $W$ is surjective, allowing us to conclude that $W$ is an adjointable
isomorphism such that $W^{-1} = W^*$.

By definition, $\D_\Om d = (\id_\sa \ot d)\D_\Om$, implying that $W(1 \ot d) = (1 \ot d) W$. Taking the adjoint of this
equation, we get $(1 \ot d^*)W^* = W^* (1 \ot d^*)$. Multiplying this equation from the left and right by $W$, it
follows that $W(1 \ot d^*) = (1 \ot d^*)W$. Applying this equation to an element of the form $1 \ot \om$, where $\om
\in \Om$, it follows that $\D_\Om d^* = (\id_\sa \ot d^*)\D_\Om$. This guarantees also the validity of the other two
commutation relations. \qed

As a consequence, $\ominv$ is invariant under the co-differential, the Dirac operator and the Laplacian.

\section{The Hodge operator}
\label{sec: Hodge operator}

{\altem To avoid tedious repetition of hypothesis, throughout this section we denote by $\omd$ a strongly
finite-dimensional $*$-differential calculus of dimension~$N$ over a Hopf $*$-algebra~$\sa$ admitting a Haar
integral~$h$. We suppose that $\dim(\ominv_N)=1$ and we denote by $\ta$ a volume element of $\omd$. We also fix a
graded, left invariant inner product $\ip{.}{.}$ on $\Om$}

\medskip
We define a linear functional $\int$ on $\Om$ by setting $\int
\om = 0$, if $\om$ is a $k$-form for which $k<N$, and $\int
\om=h(a)$, if $\om$ is an $N$-form for which $\om=a\ta$, where
$a\in \sa$. We call $\int$ the {\em integral} on $\Om$
associated to the volume element~$\ta$.

We say that $\omd$ is {\em non-degenerate} if, whenever $0\le
k\le N$ and $\om\in \Om_k$, and $\om'\om=0$, for all $\om'\in
\Om_{N-k}$, then necessarily $\om=0$.

\begin{lem} \label{thm: integral} If $\omd$ is non-degenerate,
and $\eta$ is an element of $\Om$ for which $\int \om \eta
=0$, for all $\om\in \Om$, then $\eta=0$. \end{lem}

\demo{Proof} We may clearly reduce to the case where $\eta$ is a $k$-form, for some $k\le N$. Then, if $\om$ belongs to
$\Om_{N-k}$, we have $\om \eta=a\ta$, for some element $a\in \sa$, since ${\dim(\ominv_N)=1}$. Hence, ${\int
a^*\om\eta}={\int a^*a\ta}=h(a^*a)=0$ and therefore, by faithfulness of~$h$, $a=0$. Hence, $\om\eta=0$, for all forms
$\om$ in $\Om_{N-k}$. Therefore, by non-degeneracy of $\Om$, $\eta=0$. \qed

The property enjoyed by $\int$ in this lemma is called {\em
(left) faithfulness}~\cite{KMT}.

\begin{thm} \label{thm: existence of Hodge operator} Suppose
$\omd$ is non-degenerate. Then
there exists a unique left $\sa$-linear operator $L$ on $\Om$
such that $L(\Om_k)=\Om_{N-k}$ for ${k=0,\dots,N}$ and such
that $\int \om^*L(\om')=\ip{\om'}\om$, for all $\om,\om'\in
\Om$. Moreover, $L$ is bijective. \end{thm}

\demo{Proof} Uniqueness of $L$ follows immediately from
Lemma~\ref{thm: integral}, so we confine ourselves to showing
existence. Let $\eta\in \ominv_k$ and $\xi\in \ominv_l$, where
$0\le k,l\le N$. If $k+l=N$, then $\eta^*\xi=\g(\xi,\eta)\ta$,
for a unique scalar $\g(\xi,\eta)$. Extend $\g$ to a
sesquilinear function on $\ominv$ by setting $\g(\xi,\eta)=0$
if $k+l\ne N$. Let $F$ be the unique linear operator on
$\ominv$ for which $\ip{F(\eta)}\xi=\g(\eta,\xi)$, for all
$\eta,\xi\in\ominv$. It is easily checked that
$F(\ominv_k)\subseteq \ominv_{N-k}$ and that, since $\Om$ is
non-degenerate, $F$ is injective. Since $\ominv$ is finite
dimensional, $F$ is therefore bijective. It follows that
$F(\ominv_k)= \ominv_{N-k}$. Now set $L=\id\otimes F^{-1}$,
where we identify $\Om$ with ${\sa\otimes \ominv}$, as usual.
It is clear that $L$ is a bijective left $\sa$-linear operator
on $\Om$ and that $L(\Om_k)=\Om_{N-k}$, for ${k=0,\dots,N}$.

 To show that $\int
\om^*L(\om')=\ip{\om'}\om$, for all $\om,\om'\in \Om$, we may
clearly suppose that $\om$ and $\om'$ are both $r$-forms, for
some ${r=0,\dots,N}$, using the fact that both sides of the
equation are graded sesquilinear forms. Since $\Om_r$ is a
left-covariant bi-module, and $\ominv_r$ is
finite-dimensional, we may apply Theorem~\ref{thm: basis
theorem} to choose a basis ${\ta_1,\dots,\ta_M}$ for
$\ominv_r$ and associated linear functionals $g_{jk}$ on
$\sa$, where ${j,k=1,\dots,M}$, such that
$a\ta_j={\sum_{k=1}^M \ta_k(g_{jk}*a)}$, for all $a\in \sa$.
We may also suppose that $g_{jk}(1)=\d_{jk}$. To prove $\int
\om^*L(\om')=\ip{\om'}\om$ we may make an additional
simplification and suppose now that $\om=a\ta_j$ and
$\om'=b\ta_l$, for some $a,b\in \sa$ and  some indices $j$ and
$l$  in the set ${\{1,\dots,M\}}$. Then
\begin{eqnarray*}
&\int\om^*L(\om')= \int \ta_j^*a^*bF^{-1}(\ta_l)= \int (b^*a\ta_j)^*F^{-1}(\ta_l)\\ &=\int \sum_{k=1}^M \bigl(
\ta_k(g_{jk}*(b^*a))\bigr)^*F^{-1}(\ta_l) = \int \sum_{k=1}^M (g_{jk}*(b^*a))^* \ta_k^* F^{-1}(\ta_l)\\ &=\sum_{k=1}^M
\int (\bar g_{jk}*(a^*b))\g(F^{-1}(\ta_l),\ta_k)\ta = \sum_{k=1}^M \int (\bar g_{jk}*(a^*b))\ip{\ta_l}{\ta_k}\ta\\
&=\sum_{k=1}^M h(\bar g_{jk}*(a^*b))\ip{\ta_l}{\ta_k} = \sum_{k=1}^M h(a^*b)
\bar g_{jk}(1)\ip{\ta_l}{\ta_k}\\
&=h(a^*b)\ip{\ta_l}{\ta_j}=\ip{b\ta_l}{a\ta_j}= \ip{\om'}\om.
\end{eqnarray*}
This completes the proof.\qed

The operator $L$ whose existence is shown in Theorem~\ref{thm:
existence of Hodge operator} is called the {\em Hodge
operator} on $\Om$. Let us observe that, since
$L(\ominv_k)=\ominv_{N-k}$, for $0\le k\le N$, we have
$\dim(\ominv_k)=\dim(\ominv_{N-k})$.

Observe also that $L(1)=\ta$ and $L(\ta)=1$ (these equalities
are easily verified from the definition of $F$ in the proof of
the theorem).

The following result gives a formula for the co-differential
that is familiar from classical differential geometry.

\begin{thm} \label{thm: d-L theorem} Suppose $\omd$ is non-degenerate and
$d(\ominv_{N-1})=0$. Then, if $0\le k\le N$ and $\om$ is a
$k$-form of $\Om$, we have
\[d^*\om=(-1)^kL^{-1}d L(\om).\] \end{thm}

\demo{Proof} If $k=0$, then $d^*(\om)=0$, since $d^*$ is of degree~$-1$. Also, $d L(\om)=0$, since $L(\om)$ is in
$\Om_N$ and $d(\Om_N)=0$. Hence, $d^*\om=(-1)^kL^{-1}d L(\om)$ in this case.

Suppose now that $k>0$. Clearly, $L^{-1}d L(\om)\in
\Om_{k-1}$. Hence, if ${\om'\in \Om_{k-1}}$, then
$\ip{(-1)^kL^{-1}d L(\om)}{\om'}={(-1)^k\int\om^{'*}d L(\om)}=
{\int d(\om')^*L(\om)}=\ip{\om}{d(\om')}\\=\ip{d^*\om}{\om'}$.
Here we have used the fact that $\int$ is closed; that is,
$\int d=0$. This is a consequence of our assumption that
$d(\ominv_{N-1})=0$ \cite[Corollary~4.9]{KMT}. Therefore,
$d^*\om=(-1)^kL^{-1}d L(\om)$ in this case also. \qed

\section{The Hodge decomposition}\label{sec: Hodge decomposition}

Suppose that $G=(A,\D)$ is a compact quantum group in the
sense of Woronowicz~\cite{MT, WOR2}. Let $\hat G$ denote the
set of equivalence classes of irreducible unitary
representations of $G$. For each $\a\in \hat G$, choose once
for all an irreducible unitary representation
$U^\a=(U^\a_{jk})\in {M_{N_\a}(\bC)\otimes A}$ belonging to
$\a$. Of course, $N_\a$ is the dimension of~$U^\a$. Also,
$\D(U^\a_{jk})={\sum_{r=1}^{N_\a} U^\a_{jr}\otimes
U^\a_{rk}}$, for ${j,k=1,\dots, N_\a}$.

Recall that the linear span $\sa=\sag$ of the set
${\{U^\a_{jk}}\mid {\a\in \hat G,\ {j,k=1\dots N_\a}\}}$ is a
Hopf $*$-algebra, when endowed with the restriction of the
co-multiplication $\D$ on $A$. We say that $\sag$ is the Hopf
$*$-algebra {\em associated to $G$}. Recall also that the
matrix elements~$U^\a_{jk}$ are linearly independent and
therefore form a linear basis of $\sa$.

The quantum group $G$ admits a Haar integral~$h$. When
restricted to $\sa$, $h$ is a (positive) Haar integral in the
sense of Section~\ref{sec: differential calculi}.

We denote by $\sa_\a$ the linear span of the set
${\{U^\a_{jk}}\mid {j,k=1\dots N_\a\}}$. If $\chi$ is a linear
functional on $\sa$, then $\sa_\a$ is invariant for the
operator $E_\chi={(\id\otimes \chi)\D}$. This follows from the
calculation $E_\chi(U^\a_{jk})={(\id\otimes
\chi)\D(U^\a_{jk})}= {(\id\otimes \chi)(\sum_{r=1}^{N_\a}
U^\a_{jr}\otimes U^\a_{rk}})= {\sum_{r=1}^{N_\a}
U^\a_{jr}\chi(U^\a_{rk})}$.

A {\em  $*$-differential calculus over $G$} is, by definition, a  $*$-differential calculus over the associated Hopf
$*$-algebra $\sa$. If $(\Om,d)$ is such a strongly finite dimensional calculus, set $\Om(\a)={\sa_\a\otimes \ominv}$.
Then $\Om(\a)$ is a finite-dimensional linear subspace of $\Om$ and \begin{equation} \label{eqn: alpha-decomposition}
\Om=\oplus_{\a\in\hat G} \Om(\a).\end{equation} Moreover, if we have a left invariant inner product on $\Om$, the
spaces $\Om(\a)$ are pairwise orthogonal, since the spaces $\sa_\a$ are pairwise orthogonal in
$\sa$~\cite[Theorem~7.4]{MT}.

From Equations~(\ref{eqn: tensor d-equation}) and~(\ref{eqn: tensor d*-equation})  it follows that $\Om(\a)$ is
invariant for both $d$ and~$d^*$, since $E_\chi(\sa_\a)\subseteq \sa_\a$, for every linear functional $\chi$ on $\sa$.

We shall say that an operator $T$ on a linear space $X$ is
{\em diagonalisable} if $X$ admits a linear basis consisting
of eigenvectors of $T$.

Condition~3 of the following result is our Hodge
decomposition.

\begin{thm} \label{thm: Hodge decomposition} Let $(\Om,d)$ be a strongly
finite-dimensional $*$-differential calculus over a compact quantum group $G$ and suppose we have a graded, left
invariant inner product on $\Om$ . Then

\begin{enumerate}

\item The Dirac operator $D$ and the Laplacian $\nabla$ are
diagonalisable;

\item  We have $\ker(\nabla)=\ker(D)=\ker(d)\cap\ker(d^*)$ and,
if $t$ is a positive eigenvalue of $\nabla$, then
\[\ker(\nabla-t)=\ker(D+\sqrt t)\oplus \ker(D-\sqrt t).\]

\item The space $\Om$ admits the orthogonal decomposition
\[\Om=\omhar\oplus d(\Om)\oplus d^*(\Om).\]
\end{enumerate}\end{thm}

\demo{Proof} First, let $t$ be a non-zero real number and let $t^2$ be an eigenvalue of~$\nabla$. Let $\om$ be a
corresponding eigenvector and set $\om_{\pm}=\om\pm t^{-1}D(\om)$. If $\om_+$ and $\om_-$ are non-zero, then a simple
computation shows that $\om_+$ and $\om_-$ are eigenvectors of~$D$, with $t$ and $-t$ as the corresponding eigenvalues.
If $\om_+=0$, then ${D(\om)=- t\om}$. Similarly, if $\om_-=0$, then $D(\om)=t\om$. In any case, whether $\om_+$ or
$\om_-$ are zero or non-zero, we can use the fact that ${\om=(\om_++\om_-)/2}$, to deduce that
$\ker(\nabla-t^2)=\ker(D+t)\oplus \ker(D-t)$. We shall presently  show that $\ker(\nabla)=\ker(D)$.

As we observed in the remarks preceding this theorem, each subspace $\Om(\a)$ of $\Om$ is invariant for both $d$ and
$d^*$. If $d_\a$ and $d^*_\a$ denote the restrictions of $d$ and $d^*$ to $\Om(\a)$, it is clear that $d^*_\a$ is the
adjoint of $d_\a$. Let $D_\a=d_\a+d^*_\a$ be the restriction of $D=d+d^*$ to $\Om(\a)$. Since $D_\a$ is a self adjoint
operator on a finite dimensional space its matrix with respect to a suitable basis is diagonal. This implies easily
that $\ker(D_\a^2) = \ker(D_\a)$. These observations show that $\omhar=\ker(\nabla)=\ker(D)$.  On the other hand, the
spaces $d(\Om)$ and $d^*(\Om)$ are orthogonal, since $d^2=0$. Hence, if $D(\om)=0$, then $d\om=-d^*\om$, and therefore,
$d\om=0=d^*\om$. it follows that $\ker(D)=\ker(d)\cap\ker(d^*)$. Hence, Condition~2 holds.

Also, $(d_\a d^*_\a)(d^*_\a d_\a)=0=(d^*_\a d_\a)(d_\a d^*_\a)$. Hence, $d_\a d^*_\a$ and $d^*_\a d_\a$ are commuting
(positive) normal operators on the finite-dimensional space $\Om(\a)$ and are therefore simultaneously diagonalisable.
It follows from Equation~(\ref{eqn: alpha-decomposition}) that $\Om$ admits a linear basis $(e_i)_{i\in I}$ of
simultaneous eigenvectors of $dd^*$ and $d^*d$, implying that $\nab$ is diagonalisable. Since Condition~2 holds, this
also implies that $D_\a$ is diagonalisable and we have proven Condition~1.

To prove Condition~3, we first observe that $\omhar$ is
orthogonal to the spaces $d(\Om)$ and $d^*(\Om)$. For, if
$\om\in\omhar$ and $\om',\om''\in \Om$, then
$\ip{d(\om')+d^*(\om'')}\om
=\ip{\om'}{d^*\om}+\ip{\om''}{d\om}=0$, since $d\om=d^*\om=0$,
by Condition~2.

Now let $\l_i$ and $\mu_i$ be scalars such that $(d
d^*)(e_i)=\l_i e_i$ and $(d^*d)(e_i)=\mu_i e_i$. Since
$(dd^*)(d^*d)=0$, we have $\l_i\mu_i=0$, for all $i\in I$. If
$\l_i\ne 0$, then $e_i\in d(\Om)$. Similarly, if $\mu_i\ne 0$,
then $e_i\in d^*(\Om)$. On the other hand, if $\l_i=\mu_i=0$,
then $e_i\in \omhar$. This shows that $\Om={\omhar + d(\Om)+
d^*(\Om)}$, and Condition~3 holds. \qed

\begin{cor} The orthogonal complements of
$d(\Om)$ and $d^*(\Om)$ in $\Om$ are given by
\[d(\Om)^\perp=\Om_\nabla\oplus d^*(\Om)=\ker(d^*)=\ker(dd^*)\]
and
\[d^*(\Om)^\perp=\Om_\nabla\oplus d(\Om)=\ker(d)=\ker(d^*d).\]
\end{cor}

\demo{Proof} We consider only the case of $d(\Om)^\perp$ (that
for $d^*(\Om)^\perp$ is handled the same way). It is clear
from the theorem that $d(\Om)^\perp=\Om_\nabla\oplus d^*(\Om)$
and it follows directly from the definition of $d(\Om)^\perp$
that $d(\Om)^\perp=\ker(d^*)\subseteq\ker(dd^*)$. Suppose then
$\om\in \ker(dd^*)$. By the theorem, we may write
$\om=\om_0+d\om_1+d^*\om_2$, where $\om_0\in \Om_\nabla$ and
$\om_1,\om_2\in \Om$. Then
$0=dd^*(\om)=dd^*(\om_0)+dd^*d(\om_1)=-d^*d(\om_0)+dd^*d(\om_1)$,
since $\nabla(\om_0)=0$. Hence, $d^*d(\om_0)=dd^*d(\om_1)$
belongs to $d^*(\Om)$ and to $d(\Om)$ and is therefore equal
to zero. It follows that $\nabla(d\om_1)=dd^*d(\om_1)=0$.
Consequently, $d(\om_1)=0$. Hence, $\om=\om_0+d^*(\om_2)\in
d(\Om)^\perp$. Therefore, $\ker(dd^*)=d(\Om)^\perp$. \qed

\section{Diagonalising $\nabla$ on the $0$- and $N$-forms}
\label{sec: 0-forms diagonalisation}

{\altem Throughout this section we denote by $\omd$ a strongly finite-dimensional  $*$-differential calculus of
dimension~$N$ over the Hopf $*$-algebra~$\sa$ associated to a compact quantum group~$G$. We suppose that
$\dim(\ominv_N)=1$ and that $d(\ominv_{N-1})=0$. We fix a volume element $\ta$ in $\ominv_N$ and denote by $h$ the Haar
measure of~$\sa$.We also fix a graded, left invariant inner product $\ip{.}{.}$ on $\Om$.}

\smallskip

The representations $U^\a$ are as in the first paragraph
of~Section~\ref{sec: Hodge decomposition}.

\medskip
Set $\chi=\sum_{r=1}^M\chi_r^*\chi_r$, where ${\om_1,\dots,\om_M}$ and ${\chi_1,\dots,\chi_M}$ are chosen as in
Equation~(\ref{eqn: d-chi equation}). Multiplication of linear functionals $f$ and $g$ on~$\sa$ is given, as usual, by
setting $(f g)(a)={(f\otimes g)\D(a)}$.

\begin{thm} \label{thm: Laplacian-chi theorem} If $a\in \sa$,
then $\nabla a=\chi*a$.
\end{thm}

\demo{Proof} Applying Equation~(\ref{eqn: tensor d*-equation})
to Equation~(\ref{eqn: d-chi equation}), we get
\[d^*d a=\sum_{s,r=1}^M E_{\chi_s^*\chi_r}(a)M_{\om_s}^*(\om_r)+
\sum_{r=1}^ME_{\chi_r}(a)d^*(\om_r).\] Since $d1=0$,
$d^*(\ominv_1)=0$; hence, the second sum above vanishes. It is
trivially verified that $M^*_{\om_s}(\om_r)=\delta_{rs}$.
Hence, $d^*d a=\sum_{r=1}^M E_{\chi_r^*\chi_r}(a)=E_\chi(a)$.
Since $d^*a=0$, it follows that $\nabla a=\chi*a$. \qed

The formula for $\nabla$ in the preceding theorem shows that
$\chi$ does not depend on the choice of orthonormal basis
${\om_1,\dots,\om_M}$ of~$\ominv_1$.

\begin{lem}  The operator $T^\a={(\id\otimes
\chi)\D(U^\a)\in B(H_\a)}$ is necessarily positive.
\end{lem}

\demo{Proof} If $\t$ and $\s$ are linear functionals on $\sa$,
it follows easily from the basic representation theory that
$(\id\otimes \t)(U^\a)(\id\otimes \s)(U^\a)=(\id \otimes
\t\s)(U^\a)$ and ${((\id\otimes \t)(U^\a))^*}={(\id \otimes
\t^*)(U^\a)}$. Therefore, we have $T^\a={(\id\otimes
\chi)(U^\a)}={\sum_{r=1}^M (\id\otimes
\chi_r^*\chi_r)(U^\a)}={ \sum_{r=1}^M ((\id\otimes
\chi_r)(U^\a))^*(\id\otimes \chi_r)(U^\a)}$. Thus, $T^\a$ is a
sum of positive operators and is therefore itself positive.
\qed

\begin{thm} \label{thm: 0-forms Laplacian diagonalisation}
Let $(U^\a)_\a$ be a complete family of mutually inequivalent
irreducible unitary representations of the compact quantum
group~$G$. Let $T^\a={(\id\otimes \chi)\D(U^\a)}$. Choose an
orthonormal basis $e^\a_1,\dots,e^\a_{N_\a}$ for the space
$H_\a$ on which $U^\a$ acts making $T^\a$ diagonal, so that
$T^\a e^\a_i=\l^\a_ie^\a_i$, for some elements $\l^\a_i\in
\bC$. Let $U^\a_{ij}$ be the matrix elements of $U^\a$
relative to this basis. Then $\nabla U^\a_{ij}=\l^\a_j
U^\a_{ij}$. \end{thm}

\demo{Proof} By Theorem~\ref{thm: Laplacian-chi theorem},
$\nabla U^\a_{ij}= \chi*U^\a_{ij}$. Since
$\D(U^\a_{ij})={\sum_{k=1}^{N_\a} U^\a_{ik}\otimes
U^\a_{kj}}$, we have $\chi*U^\a_{ij}={\sum_{k=1}^{N_\a}
U^\a_{ik}\chi(U^\a_{kj})}= {\sum_{k=1}^{N_\a}
U^\a_{ik}T^\a_{kj}}= {\sum_{k=1}^{N_\a}
U^\a_{ik}\l^\a_j\delta_{kj}}=\l^\a_jU^\a_{ij}$. \qed

Thus, the theorem diagonalises $\nabla_0$, the restriction of
$\nabla$ to $\sa$, since the matrix elements $U^\a_{ij}$ form
a linear basis for $\sa$. It is clear from the theorem that
the eigenvalues of $\nabla_0$ are precisely the eigenvalues of
the operators $T^\a$ and that, once we know the multiplicities
of the eigenvalues of these operators, we can calculate very
simply the multiplicities of the eigenvalues of $\nabla_0$.

\begin{cor} The eigenvalue corresponding to the
eigenvector $U^\a_{ij}$ of $\nabla$ is given by
$\chi(U^\a_{jj})=\n{dU^\a_{ij}}^2/\n{U^\a_{ij}}^2$.
\end{cor}

\demo{Proof} From the proof of the theorem we have $\nabla
U^\a_{ij}=\l^\a_j U^\a_{ij}$, where
$\l^\a_j=T^\a_{jj}=\chi(U^\a_{jj})$. Hence,
$\chi(U^\a_{jj})\ip{U^\a_{ij}}{U^\a_{ij}}= \ip{\nabla
U^\a_{ij}}{U^\a_{ij}}=\ip{d^*d U^\a_{ij}}{U^\a_{ij}}=
\ip{dU^\a_{ij}}{dU^\a_{ij}}$. \qed

Woronowicz defines the compact quantum group~$G$ to be {\em
connected} if Condition~1 of the following theorem is
satisfied. The theorem provides two other equivalent
conditions.

\begin{thm}The following are equivalent
conditions:

\begin{enumerate}
\item For any element $a\in \sa$, $d a=0$ if, and only if,
$a=\l 1$, for some $\l\in \bC$;

\item For any element $a\in \sa$, $\nabla a=0$ if, and only
if, $a=\l 1$, for some $\l\in \bC$;

\item If $U^\a\ne U^{\a_0}$, the trivial unitary representation of~$G$,
then the operator $T^\a={(\id\otimes\chi)\D(U^\a)}$ is
invertible.
\end{enumerate}\end{thm}

\demo{Proof} Condition~2 clearly implies~1. Conversely, if
Condition~1 holds, then so does~2, since
$\ker(\nabla)=\ker(d)\cap\ker(d^*)$ and $d^*a=0$, for all
$a\in \sa$.

Suppose now that Condition~2 holds and suppose that $U^\a\ne
U^{\a_0}$ and that $T^\a$ is not invertible. Then $\l^\a_j=0$,
for some~$j$ and therefore $\nabla U^\a_{jj}=\l^\a_j
U^\a_{jj}=0$. Hence, by Condition~2, there exists a scalar $c$
for which $U^\a_{jj}=c1$. However, since $U^\a\ne U^{\a_0}$,
the element $U^\a_{jj}$ is orthogonal to~$1$ in $\sa$.
Therefore, $c=h(1^*U^\a_{jj})=0$. Thus, $U^\a_{jj}=0$, an
impossibility. To avoid contradiction it must therefore be the
case that $T^\a$ is invertible. Hence, Condition~2 implies~3.

Finally, suppose that Condition~3 holds. Let $a\in \sa$ be
such that $\nabla a=0$. Because the $U^\a_{ij}$ form a linear
basis of $\sa$, we may write ${a=\sum_\a\sum_{i,j=1}^{N_\a}
c^\a_{ij} U^\a_{ij}}$, for some scalars~$c^\a_{ij}$. Since
$0=\nabla a= {\sum_\a\sum_{i,j=1}^{N_\a} c^\a_{ij} \l^\a_j
U^\a_{ij}}$, we have $c^\a_{ij} \l^\a_j=0$, for all $\a$, $i$
and~$j$. Hence, if $U^\a\ne U^{\a_0}$, then $\l^\a_j\ne 0$,
and therefore $c^\a_{ij}=0$. It follows that $a=c1$, for
$c=c^{\a_0}_{11}$. Hence, Condition~3 implies~2 \qed

We turn now to finding the eigenvalues and eigenvectors of the
Laplacian on the highest dimensional forms, the $N$-forms. The
situation here turns out to be analogous to that for the
lowest-dimensional forms, the $0$-forms.

We first need to introduce a new linear functional $\g$ on
$\sa$, an analogue for the case of the $N$-forms of the
functional $\chi$ used for the case of the $0$-forms. It is
clear, since $\ominv_N=\bC \ta$, that
$M_{\om_r}M^*_{\om_s}(\ta)=c_{rs}\ta$, for some scalars
$c_{rs}$. Since $c_{rs}=\ip{M_{\om_r}M^*_{\om_s}\ta}\ta=
\ip{M^*_{\om_s}\ta}{M^*_{\om_r}\ta}$, the matrix $(c_{rs})$ is
positive. We set $\g=\sum_{r,s=1}^M c_{rs}\chi_r\chi^*_s$.

\begin{thm} \label{thm: Laplacian-gamma theorem} If $a\in \sa$, then
$\nabla(a\ta)=(\g*a)\ta$.
\end{thm}

\demo{Proof} Note first that $d(a\ta)=0$, since $d(\Om_N)=0$.
Hence, $\nabla(a\ta)=dd^*(a\ta)$. Using Equation~(\ref{eqn:
tensor d*-equation}), we get
$d^*(a\ta)={ad^*(\ta)+\sum_{s=1}^ME_{\chi_s^*}(a)M^*_{\om_s}(\ta)}$.
The first term here is equal to zero, since $d^*(\ta)=0$; this
is an immediate consequence of the fact that
$d(\ominv_{N-1})=0$. It follows, using Equation~(\ref{eqn:
tensor d-equation}), that $dd^*(a\ta)={\sum_{s=1}^M
E_{\chi_s^*}(a)dM^*_{\om_s}(\ta)}+{\sum_{r,s=1}^M
E_{\chi_r\chi_s^*}(a)M_{\om_r}M^*_{\om_s}(\ta)}$. Again using
the fact that $d(\ominv_{N-1})=0$, we see that the first sum
vanishes. Since $M_{\om_r}M^*_{\om_s}(\ta)=c_{rs}\ta$, we have
$\nabla(a\ta)=dd^*(a\ta)={\sum_{r,s=1}^M
c_{rs}E_{\chi_r\chi_s^*}(a)\ta}=(\g*a)\ta$. \qed

The formula for $\nabla$ restricted to the $N$-forms in the
preceding proposition shows that $\g$ does not depend on the
choice of orthonormal basis ${\om_1,\dots,\om_M}$
of~$\ominv_1$.

\begin{lem}   The operator $S^\a=(\id\otimes \g)(U^\a)\in B(H_\a)$
is necessarily positive.
\end{lem}

\demo{Proof} We have $S^\a= {\sum_{r,s=1}^M c_{rs}(\id\otimes
\chi_r)(U^\a)((\id\otimes \chi_s)(U^\a))^*}$. Also, since the
matrix $(c_{rs})$ is positive, it can be written as a sum of
matrices ${c^1,\dots,c^K}$, with each such matrix is of the
form $c^k_{rs}= c^k_r\bar c^k_s$, for some scalars $c^k_r$.
Hence,
\begin{eqnarray*}
&S^\a= {\sum_{k=1}^K\sum_{r,s=1}^M c^k_r\bar c^k_s(\id\otimes
\chi_r)(U^\a) ((\id\otimes \chi_s)(U^\a))^*}\\ &={\sum_{k=1}^K
(\sum_{r=1} c^k_r (\id\otimes \chi_r)(U^\a))(\sum_{r=1} c^k_r
(\id\otimes \chi_r)(U^\a))^*}.
\end{eqnarray*}
Thus, $S^\a$ is a sum of positive operators and is therefore
itself positive.~\qed

\begin{thm} \label{thm: N-forms Laplacian diagonalisation}
Let $(U^\a)_\a$ be a complete family of mutually inequivalent
irreducible unitary representations of the compact quantum
group~$G$. Set $S^\a=(\id\otimes \g)(U^\a)$ and choose an
orthonormal basis $e^\a_1,\dots,e^\a_{N_\a}$ for the space
$H_\a$ on which $U^\a$ acts making $S^\a$ diagonal, so that
$S^\a e^\a_i=\l^\a_ie^\a_i$, for some elements $\l^\a_i\in
\bC$. Let $U^\a_{ij}$ be the matrix elements of $U^\a$
relative to this basis. Then $\nabla(U^\a_{ij}\ta)=\l^\a_j
U^\a_{ij}\ta$.
\end{thm}

\demo{Proof} By Theorem~\ref{thm: Laplacian-gamma theorem},
$\nabla(U^\a_{ij}\ta)= (\g*U^\a_{ij})\ta$. Since
$\g*U^\a_{ij}={\sum_{k=1}^{N_\a} U^\a_{ik}\g(U^\a_{kj})}\\=
{\sum_{k=1}^{N_\a} U^\a_{ik}S^\a_{kj}}= {\sum_{k=1}^{N_\a}
U^\a_{ik}\l^\a_j\delta_{kj}}=\l^\a_jU^\a_{ij}$, we have
$\nabla(U^\a_{ij}\ta)=\l^\a_j U^\a_{ij}\ta$, as required. \qed

Remarks analogous to those made after Theorem~\ref{thm:
0-forms Laplacian diagonalisation} apply to Theorem~\ref{thm:
N-forms Laplacian diagonalisation}.

We show now that a formula holds for the Laplacian applied to
general $k$-forms that is analogous to those we obtained in
Theorems~\ref{thm: Laplacian-chi theorem} and~\ref{thm:
Laplacian-gamma theorem} for 0-forms and $N$-forms,
respectively.

Before proceeding, note that $\sa_\a$ is finite dimensional,
of dimension $N_\a^2$. Let $k$ be an integer such that $0\le
k\le N$. The linear span $\Om_{\a,k}$ of the forms $a\om$,
where $a\in \sa_\a$ and $\om \in \ominv_k$, is linearly
isomorphic to the tensor product ${\sa_\a\otimes \ominv_k}$,
where the isomorphism maps $a\om$ onto $a\otimes \om$. We
shall use this isomorphism to identify these two spaces. Using
the fact that the matrix entries of $U^\a$ are orthogonal to
the matrix entries of~$U^\b$, if $\a\ne \b$, it is clear that
$\Om_k=\oplus_\a \Om_{\a,k}$, as an orthogonal sum.

Equations~(\ref{eqn: tensor d-equation}) and (\ref{eqn: tensor
d*-equation}) imply that the space $\Om_k={\sa\otimes
\ominv_k}$ is invariant for $\nabla$; we may therefore write
$\nabla ={\nabla_0\oplus\cdots\oplus\nabla_N}$, where
$\nabla_k$ is the restriction of $\nabla$ to~$\Om_k$.
Similarly, the spaces $\Om_{\a,k}$ are also invariant for
$\nabla$. Hence, $\nabla_k=\oplus_\a \nabla_{\a,k}$, where
$\nabla_{\a,k}$ is the restriction of $\nabla$ to
$\Om_{\a,k}$. We therefore reduce the eigenvalue problem for
$\nabla$ to one for the finite-dimensional operator
$\nabla_{\a,k}$.

Recall, for the proof of the following theorem, the elementary
fact that if $(e_i)_{i\in I}$ is a finite linear basis for a
linear space $X$, the corresponding system of matrix units
$(e_{ij})_{i,j\in I}$ in $B(X)$ is defined by setting
$e_{ij}(e_k)=\d_{jk}e_i$ and that if $T\in B(X)$ has matrix
$(\l_{ij})$ relative to the basis $(e_i)$, then
$T=\sum_{i,j\in I} \l_{ij}e_{ij}$.

\begin{thm} \label{thm: block matrix theorem} Suppose $k$ is an integer
such that ${0\le k\le N}$ and let ${\eta_1,\dots,\eta_L}$ be a
linear basis for $\ominv_k$.

\begin{enumerate}

\item For $p,q=1,\dots,L$, there exist unique linear functionals
$g_{pq}$ on
$\sa$ such that
\[ \nabla(a\eta_p)=\sum_{q=1}^L (g_{qp}*a )\eta_q ,\]
for all $a\in \sa$.

\item If $\a\in \hat G$, the restriction $\nabla_{\a,k}$ of
$\nabla$ to $\Om_{\a,k}$ is similar to a direct sum of $N_\a$
copies of the block matrix
\[
\left[
\begin{array}{cccc}
       g_{11}(U^\a)&g_{12}(U^\a)&\dots&g_{1L}(U^\a)\\
       g_{21}(U^\a)&g_{22}(U^\a)&\dots&g_{2L}(U^\a)\\
       \vdots&\vdots&\ddots&\vdots\\
       g_{L1}(U^\a)&g_{L2}(U^\a)&\dots&g_{LL}(U^\a)
\end{array}\right],
\]
where $g_{pq}(U^\a)$ is the $N_\a\times N_\a$-matrix with $i
j$th entry $g_{pq}(U^\a_{ij})$.

\end{enumerate}

\end{thm}

\demo{Proof} Let $M^k_r$ denote the operator from $\ominv_k$
to $\ominv_{k+1}$ of left multiplication by $\om_r$. Let
$\chi_0$ be the co-unit of $\sa$ and let $M^k_0$ be the
restriction of $d$ to $\ominv_k$. Finally, set
$E_r=E_{\chi_r}$. By Theorem~\ref{thm: d-adjointability
theorem}, $\nabla_k={\sum_{r,s=0}^M E_rE^*_s\otimes
M^{k-1}_r(M^{k-1}_s)^*}+ {\sum_{r,s=0}^M E^*_rE_s\otimes
(M^k_r)^*M^k_s}$. Since the operators $M^{k-1}_r(M^{k-1}_s)^*$
and $(M^k_r)^*M^k_s$ are linear combinations of the matrix
units $\eta_{pq}$ corresponding to the basis ${\eta_1,
\dots,\eta_L}$ of $\ominv_k$, it follows that
$\nabla_k={\sum_{p,q=1}^L E_{g_{pq}}\otimes \eta_{pq}}$, for
linear functionals $g_{pq}$ that are linear combinations of
the functionals $\chi_r\chi^*_s$ and $\chi^*_r\chi_s$. This
proves Condition~1 (uniqueness of the functionals $g_{pq}$ is
clear).

In the remainder of this proof, in summations the indices $p,
q$ will range over the numbers ${1,\dots,L}$ and the indices
$i,j$ over the range ${1,\dots,N_\a}$. For ${m\in
\{1,\dots,N_\a\}}$, we have
$\nabla(U^\a_{mi}\eta_p)={\sum_q(g_{qp}*U^\a_{mi})\eta_q}={\sum_{q,j}
U^\a_{mj}g_{qp}(U^\a_{ji})\eta_q}$. Hence, if $Y=\Om_{\a,k}$
and $X_m$ is the linear span of the elements
${U^\a_{m1},\dots,U^\a_{mN_\a}}$, then $Y_m={X_m\otimes
\ominv_k}$ is an invariant subspace for $\nabla$ and
$Y={Y_1\oplus\cdots\oplus Y_{N_\a}}$. We write $T_m$ for the
restriction of $\nabla$ to $Y_m$. Then
$\nabla_{\a,k}={T_1\oplus\cdots\oplus T_{N_\a}}$. We shall
show that each summand is similar to the block matrix
displayed in Condition~2 and this will complete the proof of
the theorem.

Let $e_i=U^\a_{mi}$, for ${i=1,\dots,N_\a}$. Then the $e_i$
form a basis for $X_m$ and therefore the elements ${e_i\otimes
\eta_p}$ form a basis for $Y_m$. Let $e_{ij}$, $\eta_{pq}$ and
$f_{(i,p),(j,q)}$ form a system of matrix units corresponding
to the bases $e_i$, $\eta_p$ and ${e_i\otimes \eta_p}$,
respectively. It is elementary to verify that
$f_{(i,p),(j,q)}={e_{ij}\otimes \eta_{pq}}$. Since
$\nabla(U^\a_{mj}\eta_q)={\sum_{p,i}
U^\a_{mi}g_{pq}(U^\a_{ij})\eta_p}$, we have
\[ T_m= \sum_{ijpq}
g_{pq}(U^\a_{ij})f_{(i,p),(j,q)}=
\sum_{pq}(\sum_{ij}g_{pq}(U^\a_{ij})e_{ij})\otimes
\eta_{pq}.\] Hence, $T_m$ has matrix with respect to the basis
${e_i\otimes \eta_p}$ of $Y_m$ the block matrix displayed in
Condition~2 and therefore $T_m$ is similar to that matrix, as
required.~\qed

Thus, the eigenvalue problem for the Laplacian $\nabla$ has
been reduced to the finite-dimensional problem of finding the
eigenvalues of the block matrix in Condition~2 of
Theorem~\ref{thm: block matrix theorem}. However, this matrix
is of order $N_\a L$; it may be difficult to work with in
practice, since, in general, $N_\a$ can be arbitrarily large.
In Section~\ref{sec: eigenvalues of nabla} we present an
alternative  approach in a special case.

\section{The inner product on $\ominv$}
\label{sec: inner product}

{\altem In this section we make the same assumptions and notational conventions as in the first paragraph of
Section~\ref{sec: Hodge operator} without fixing an inner product on $\Om$. In addition, we suppose that $\Om$ is
non-degenerate and that $\int$ is the integral associated to~$\ta$.}

\medskip
In general, the Laplacian $\nabla$ and the Hodge operator $L$
need not commute, as one can see by looking at examples.
However, when they do commute, the problem of computing the
eigenvalues of $\nabla$ is effectively halved. For, if $\om\in
\Om_k$ and $\nabla \om =\l \om$, for some scalar $\l$, then
$\nabla L(\om)=\l L(\om)$ also. Thus, if one computes the
eigenvalues of $\nabla$ on $\Om_k$, one knows them for
$\nabla$ on $\Om_{N-k}$ also, where $N$ is the dimension
of~$\Om$.

Up to this point, we have not specified how the inner product on $\ominv$ should be chosen, except to say that it
should be graded. It turns out that one can indeed choose the inner product on $\ominv$ in such a way as to ensure that
$\nabla L=L\nabla$. We turn now to showing this.

Choose an inner product $\ip \cdot\cdot'$ on $\ominv$ and let
$L'$ be the corresponding Hodge operator on~$\Om$. We obtain a
new graded inner product $\ip{\cdot}\cdot$ on $\ominv$, by
defining $\ip{\om_1}{\om_2}$, for $\om_1,\om_2\in \ominv_k$,
as follows:

\begin{enumerate}

\item If $k<N/2$, set
$\ip{\om_1}{\om_2}=\ip{\om_1}{\om_2}'$;

\item If $k>N/2$, set
$\ip{\om_1}{\om_2}=\ip{(L')^{-1}(\om_1)}{(L')^{-1}(\om_2)}'$;

\end{enumerate}

For the moment, we set aside the case of $k=N/2$, which
occurs, of course, only if $N$ is even; we suppose only that
the new inner product has been specified in some way on
$\ominv_{N/2}$ in this case.

Notice that the new inner product on $\ominv_N$ is the same as
the old (since $L'(1)=\ta$). Hence, $\ta$ is still a volume
element for $\Om$ with its new inner product.

Let $L$ be the Hodge operator on $\Om$ associated to this new
inner product.

Observe that, if $k<N/2$, then $L=L'$ on $\Om_k$. This follows
immediately from faithfulness of~$\int$.

On the other hand, if $k>N/2$, then ${L=(-1)^{k(N-k)}(L')^{-1}}$ on $\Om_k$. To see this, we need only show that
$L(\om)={(-1)^{k(N-k)}(L')^{-1}(\om)}$, for $\om\in \ominv_k$. Hence, we need only show ${\int \eta^*
(L(\om)-(-1)^{k(N-k)}(L')^{-1}(\om))=0}$, for all $\eta\in \ominv_k$, by faithfulness of~$\int$. To see this equality,
observe first that if $\om\in \ominv_{N-k}$ and $\eta\in \ominv_k$, then
\begin{eqnarray*}
&\ip{(L')^{-1}(\om)}{\eta}'={\int \eta^*\om}=\overline{\int (\eta^*\om)^*}=(-1)^{k(N-k)}\overline{\int \om^*\eta}=\\ &
{(-1)^{k(N-k)}\overline{\ip{(L')^{-1}(\eta)}{\om}'}} ={(-1)^{k(N-k)}\ip{\om}{(L')^{-1}(\eta)}'}. \end{eqnarray*} Hence,
if $\eta,\om\in \ominv_k$, then ${\int \eta^* L(\om)}=\ip{\om}\eta=\ip{(L')^{-1}(\om)}{(L')^{-1}(\eta)}'
=(-1)^{k(N-k)}\ip{(L')^{-2}(\om)}{\eta}'= {(-1)^{k(N-k)}\int \eta^*(L')^{-1}(\om)}$. Therefore, as required, we have
${\int \eta^* (L(\om)-(-1)^{k(N-k)}(L')^{-1}(\om))=0}$.

It follows from these calculations that
$L^2(\om)=(-1)^{k(N-k)}\om$, for all $\om\in \Om_k$, if $k\ne
N/2$. We show now how to choose the inner product on
$\ominv_k$, where $k=N/2$, in such a way that, in this case
also, $L^2(\om)=(-1)^{k(N-k)}\om$, for all $\om\in \Om_k$:

The sesquilinear form, ${(\eta,\om)\mapsto \int\om^*\eta}$, is
non-degenerate on $\ominv_k$, by faithfulness of~$\int$. Also,
this form is Hermitian if $(-1)^{k(N-k)}=1$ and anti-Hermitian
if $(-1)^{k(N-k)}=-1$. Set $\e_k=1$ in the first case and
$\e_k=i$ in the second. Then there exists a basis
${e_1,\dots,e_m}$ for $\ominv_k$ such that $\int
e_q^*e_p=\pm\e_k\d_{pq}$, for all ${p,q=1,\dots,m}$. Now
choose an inner product on $\ominv_k$ making ${e_1,\dots,e_m}$
orthonormal. For this inner product, the operator $L$ on
$\ominv_k$ defined by mapping $e_p$ onto $(-1)^{k(N-k)}(\int
e_p^*e_p)e_p$, for ${p=1,\dots,m}$, satisfies the condition
${\int \eta^*L(\om)}=\ip\om\eta$; hence, $L$ is the Hodge
operator on $\ominv_k$ and extends in an obvious way to the
Hodge operator on~$\Om_k$. Clearly,
$L^2(\om)=(-1)^{k(N-k)}\om$, for all $\om\in\Om_k$.

It is clear from our construction of $L$ that it is a unitary operator on~$\Om$.

 We sum up what we have shown here in the
following theorem and draw out some useful consequences.

\begin{thm} \label{thm: L-properties theorem} Let $\omd$ be a
non-degenerate strongly finite-dimensional $*$-differential
calculus of dimension~$N$ over a Hopf algebra~$\sa$ admitting
a Haar integral. Suppose that ${\dim(\ominv_N)=1}$ and
$d{(\ominv_{N-1})=0}$. Then $\ominv$ admits an inner product
such that the corresponding Hodge operator~$L$ is a unitary
and satisfies the following conditions

\begin{enumerate}

\item $L^2(\om)=(-1)^{k(N-k)}\om$, for all $\om\in \Om_k$;

\item $L d(\om)=(-1)^{k+1}d^*L(\om)$, for all $\om\in \Om_k$;

\item $d^*dL=Ldd^*$, $dd^*L=Ld^*d$ and $L\nabla=\nabla L$.
\end{enumerate}
\end{thm}

\demo{Proof} Condition 1 was proved above. To show~2,  first
observe that we may suppose that $k<N$. Now observe that, by
Theorem~\ref{thm: d-L theorem},
$d\om=(-1)^{N-k}Ld^*L^{-1}(\om)$. Hence,
\begin{equation}\label{eqn: d-L equation}
L^{-1}d\om=(-1)^{N-k}d^*L^{-1}(\om).
\end{equation}
Now, by Condition~1, we have
\begin{eqnarray*}
&{(-1)^{(k+1)(N-k-1)}Ld\om}= {L^{-1}d\om} ={(-1)^{N-k}d^*L^{-1}(\om)}\mbox{ (by Equation~(\ref{eqn: d-L equation}))}\\
&=(-1)^{N-k}(-1)^{k(N-k)}d^*L(\om)\mbox{ (by Condition~1 again).}
\end{eqnarray*}
Hence,
$Ld\om=(-1)^{(k+1)(N-k-1)}(-1)^{N-k}(-1)^{k(N-k)}d^*L(\om)
=(-1)^{k+1}d^*L(\om)$. Thus, Condition~2 holds.

Again with $\om\in \Om_k$, we get $d^*d
L(\om)={(-1)^kd^*Ld^*\om} =(-1)^k(-1)^kLdd^*\om$, by
Theorem~\ref{thm: d-L theorem} and Condition~2. Hence,
$d^*dL=Ldd^*$. It follows from Condition~1 that $dd^*L=Ld^*d$.
Hence, $\nabla L=dd^*L+d^*d L= Ld^*d+Ldd^*=L\nabla$, so
Condition~3 holds. \qed

\section{One-parameter representations induced by \newline
twist automorphisms} \label{sec: one-parameter reps}

{\altem Throughout this section the same standing assumptions
and notational conventions will be used as specified in the
first paragraph of Section~\ref{sec: 0-forms diagonalisation},
with the added assumption that $\Om$ is non-degenerate. We
will also make use of the notations introduced at the
beginning of Section~\ref{sec: Hodge decomposition}.

The sets of nonnegative and positive integers will be denoted
by $\N$ and $\N_0$, respectively.}

\medskip
In this section we shall need a number of results from the
general theory of exponential one-parameter groups of
automorphisms. Since these are not completely standard, we
have gathered them together in an appendix to which we shall
make frequent reference. We shall therefore use in this
section the terminology and notations of the appendix.

\medskip

We begin by introducing the two most important one-parameter
groups connected to the Hopf algebra~$\sa$.
By~\cite[Theorem~5.6]{WOR2} there exist unique exponential
one-parameter groups~$\r$ and~$\t$ on~$\sa$ such that $h(a
b)=h(\r_i(b)a)$, for all $a,b\in \sa$ and $\t_i=\k^{-2}$. We
call $\r$ the {\em modular group} of $\sa$ and of $h$, and
$\t$ the {\em scaling group} of $\sa$.

There is a family of linear functionals $(f_z)_{z\in\C}$
on~$\sa$ such that
\[\rho_z(a) =
(f_{iz} \ot \id_{\sa} \ot f_{iz})\D^{(2)}(a) \hspace{0.5cm} \text{and} \hspace{0.5cm} \tau_z(a) = (f_{iz} \ot \id_{\sa}
\ot f_{-iz})\D^{(2)}(a),\] for all $z \in \C$ and $a \in \sa$. Here we use the usual convention that $\D^{(2)} = (\D
\ot \id_{\sa})\D = (\id_{\sa} \ot \D)\D$. (See \cite[Theorem~5.6]{WOR2} for these functionals). It follows that
$\r_z(\sa_\a)=\sa_\a$ and $\t_z(\sa_\a)=\sa_\a$, for all $\a\in\hat G$.

We note the following facts: $h \, \rho_z = h \, \tau_z = h$,
for all $z \in \C$. Also, $(\tau_z \ot \rho_z)\D = (\rho_z \ot
\tau_{-z})\D = \D\,\rho_z$ and $(\tau_z \ot \tau_z)\D =
(\rho_z \ot \rho_{-z})\D = \D \, \tau_z$.

Recall that there exists  a unique $*$-anti-automorphism $R$ on $\sa$ such that $\k = R \, \tau_{-i/2} = \tau_{-i/2} \,
R$ and $R^2 = \id_{\sa}$. This decomposition of $\k$ is called the {\em polar decomposition} of $\k$; $R$ is called the
{\em unitary co-inverse} of $\sa$.

Note that $\flip(R \ot R)\D = \D R$, where $\flip$ is the flip map. Also, $h\, R = h$. We refer the reader to section~7
of \cite{WOR4} for these results.


\begin{lem} \label{one.lem2}
Suppose that $\pi : \sa \rightarrow \sa$ is a surjective linear map. If $\pi$ is multiplicative (respectively,
anti-multiplicative) and $h \, \pi = h$, then $\pi \, \rho_z = \rho_z \, \pi$ (respectively, $\pi \, \rho_z = \rho_{-z}
\, \pi$), for all $z \in \C$.
\end{lem}

\demo{Proof} We prove only the multiplicative case---the
anti-multiplicative case is proved similarly. For $a,b \in
\sa$, we have $h\bigl(\rho_i(\pi(b)) \pi(a)\bigr) = h(\pi(a)
\pi(b)) = h(a b) = h( \rho_i(b) a) = h\bigl(\pi(\rho_i(b))
\pi(a)\bigr)$. Therefore, surjectivity of $\pi$ and
faithfulness of $h$ imply that $\rho_i(\pi(b)) =
\pi(\rho_i(b))$. Hence, Lemma~\ref{one.lem1} implies that
$\rho_z \, \pi = \pi \, \rho_z $, for all $z \in \C$. \qed

The following special cases should be noted: $\tau_z \, \rho_y
= \rho_y \, \tau_z$,  $\rho_y \, R = R \, \rho_{-y}$ and
$\rho_y \, \k = \k \,  \rho_{-y}$, for all $y,z \in \C$. Note
also that ${\tau_z \, R = R \, \tau_z}$ and  $\tau_z \, \k =
\k \, \tau_z$.

\begin{lem}
Suppose that $\phi : \sa \rightarrow \sa$ is a surjective
algebra homomorphism  such that $h \phi = h$ and such that
$h(a^* \phi(a))$ is real, for all $a \in A$. Then
$\phi(\phi(a^*)^*) = a$, for all $a \in \sa$.
\end{lem}

\demo{Proof} The function $\sa \times \sa \rightarrow \C :
(a,b) \mapsto h(b^* \phi(a))$ is a hermitian form on $\sa$.
Hence, for $a,b \in \sa$, $h(\phi(\phi(a^*)^*)\, \phi(b)) = h
(\phi(a^*)^* b) = {h(b^* \phi(a^*))^-} =  h(a \phi(b))$.
Therefore, by faithfulness of $h$, $\phi(\phi(a^*)^*) = a$.
\qed

Let $\int$ denote the integral on $\Om$ associated to $\ta$. It is a twisted graded trace on $(\Om,d)$,
by~\cite[Corollary~4.9]{KMT}. This means there is an algebra automorphism $\s$ of~$\Om$ of degree zero such $\s d=d\s$
and ${\int \om'\om}=(-1)^{kl}{\int \s(\om)\om'}$, for all nonnegative integers $k$ and $l$ and all $\om\in \Om_k$ and
$\om'\in \Om_l$. Faithfulness of $\int$ implies that $\s$ is unique. We call $\s$ the {\em twist automorphism} of
$\int$. Note that $\int\s=\int$. Also, $\s(\s(\om^*)^*) = \om$, for all $\om \in \Om$. In order to prove this we may
and will assume that $\om \in \Om_k$. Since $\int$ is self adjoint, we get for all $\eta \in \Om_{N-k}$,
\begin{eqnarray*}
& & \int \eta \, \s(\s(\om)^*)^* = (-1)^{k(N-k)}\, \bigl(\,\int \s(\s(\om)^*)
\,\eta^*\,\bigr)\bar{\rule[.5ex]{0ex}{1.3ex}}
\\ & & \hspace{1ex} = \bigl(\int \eta^* \s(\om)^*\,\bigr)\bar{\rule[.5ex]{0ex}{1.3ex}} =
(-1)^{k(N-k)}\, \int \s(\om) \, \eta = \int \eta \, \om .
\end{eqnarray*}
Therefore the faithfulness of $\int$ implies that $\s(\s(\om^*)^*) = \om$. The restriction of $\s$ to $\sa$ will be
denoted by $\s_{\!\sa}$; it is an automorphism on~$\sa$.

The proof of the following lemma is modeled on that of
\cite[Proposition~3.14]{VD}.

\begin{lem} \label{gjm: kappa lemma}
$(\k^{-2} \ot \s)\D_\Om = \D_\Om \, \s$ and
$\s(\ominv)=\ominv$.
\end{lem}

\demo{Proof} Let $\om \in \Om_k$ and $\eta \in \Om_{N-k}$. By \cite[Theorem~4.4]{KMT}, we have
\begin{eqnarray*}
& & \hspace{-3ex} (\id_{\sa} \ot \int)\bigl((\k^{-2} \ot
\s)(\D_\Om(\om)) (1 \ot \eta)\bigr)
\\ & & \hspace{1ex} = \k^{-2}\bigl(\,(\id_{\sa}
\ot \int)\bigl((\id_{\sa} \ot \s)(\D_\Om(\om)) (1 \ot
\eta)\bigr)\,\bigr)
\\ & & \hspace{1ex} = (-1)^{k(N-k)}\,\k^{-2}\bigl(\,(\id_{\sa} \ot \int)\bigl( (1 \ot \eta)
\D_\Om(\om)\bigr)\,\bigr)
\\ & & \hspace{1ex} = (-1)^{k(N-k)}\,\k^{-1}\bigl(\,(\id_{\sa} \ot \int)\bigl(\D_\Om(\eta) (1 \ot
\om)\bigr)\,\bigr)
\\ & & \hspace{1ex} =  \k^{-1}\bigl(\,(\id_{\sa} \ot \int)\bigl((1 \ot
\s(\om))\D_\Om(\eta)\bigr)\,\bigr)
 = (\id_{\sa} \ot \int)\bigl(\D_\Om(\s(\om)) (1 \ot \eta)\bigr) \ ,
\end{eqnarray*}
Consequently, faithfulness of $\int$ implies that $(\k^{-2} \ot \s)\D_\Om(\om) = \D_\Om(\s(\om))$. This equation
clearly implies in turn that $\s(\Om^{\text{inv}}) \subseteq \Om^{\text{inv}}$. Since $\s(\s(\om^*)^*) = \om$, we have
$\s(\ominv)=\ominv$. \qed

\begin{thm} For all $z \in
\C$, $h\,\s_{\!\sa} = h$, $\s_{\!\sa} \, \rho_z = \rho_z \, \s_{\!\sa}$ and $\s_{\!\sa} \, \tau_z  = \tau_z \,
\s_{\!\sa}$.
\end{thm}

\demo{Proof} Theorem~4.5 of \cite{KMT} implies that $h(\s(a)) = \int \s(a) \ta = \int \ta a = h(a)$, for every $a \in
\sa$, proving the first equation holds. The second equation now follows from Lemma~\ref{one.lem2}.

Using Lemma~\ref{gjm: kappa lemma}, we see that, for $z\in
\C$,
\begin{eqnarray*}
& & (\rho_z\, \k^{-2} \ot \tau_{-z}\, \s_{\!\sa}) \D = (\rho_z \ot \tau_{-z})\D \s_{\!\sa} = \D \rho_z \s_{\!\sa} = \D
\s_{\!\sa} \rho_z
\\ & & \spat = (\k^{-2} \ot \s_{\!\sa})\D \rho_z = (\k^{-2}
\,\rho_z \ot \s_{\!\sa}\, \tau_{-z})\D = (\rho_z \, \k^{-2} \ot \s_{\!\sa} \, \tau_{-z} )\D \ ,
\end{eqnarray*}
implying that $(\id_{\sa} \ot \tau_{-z}\, \s_{\!\sa})\D = (\id_{\sa} \ot \s_{\!\sa}\,\tau_{-z})\D$. Since $\sa$ is the
linear span of elements of the form $(f\ot \id_{\sa})\D(a)$, where $a \in \sa$ and $f \in \sa'$, we conclude that
$\tau_{-z} \, \s_{\!\sa} = \s_{\!\sa} \, \tau_{-z}$. \qed

We extend the concept of a one-parameter group from a
$*$-algebra to our differential calculus~$\Om$. A map $\b$
from $\C$ into the space of algebra automorphisms  on $\Om$ is
a {\em one-parameter group of automorphisms} if
$\b_{y+z}=\b_y\b_z$ and $\b(\om)^*=\b_{\bar z}(\om^*)$, for
all $y,z\in\C$ and $\om\in \Om$; furthermore, it is {\em
exponential} if the map, ${\C\to \Om : z\mapsto \b_z(\om)}$,
is of finite exponential type, for all $\om\in \Om$.

We shall call a family of elements $(e_j)_{j \in J}$ in $\sa$
{\em homogeneous} if each of the elements $e_j$ belongs to
some space $\sa_\a$, for some $\a\in \hat G$.

\begin{thm} \label{one.prop2}
Suppose that $h(a^* \s(a)) \geq 0$, for all $a \in \sa$. Then there exists a unique exponential one-parameter group
$\b$ on $\Om$ such that $\b_i = \s$. Moreover, there is a homogeneous orthonormal basis $(e_j)_{j \in J}$ for $\sa$ and
there are positive numbers $(\lambda_j)_{j \in J}$ such that $\b_z(e_j) = \lambda_j^{-iz} \, e_j$, for all $j \in J$
and $z \in \C$.

Also, $\b_z$ is of degree zero,  $d\, \b_z = \b_z\, d$ and
$\int \b_z = \int$. Finally, $(\tau_z \ot \b_z)\D_\Om = \D_\Om
\b_z$ and $\b_z(\ominv) = \ominv$.
\end{thm}

\demo{Proof} Uniqueness of $\b$ is obvious, so we concern
ourselves only with existence. If $\a \in \hat{G}$,
$\k^{-2}(\sa_\a)=\t_i(\sa_\a)=\sa_\a$. Using Lemma~\ref{gjm:
kappa lemma}, we get,  for $i,j= 1,\ldots,N_\a$,
$$\D(\s(U^\a_{ij})) = (\k^{-2} \ot \s_{\!\sa})\D(U^\a_{ij}) =
\sum_{k=1}^{N_\a} \k^{-2}(U^\a_{ik}) \ot \s(U^\a_{kj})\ .$$ Applying $\id_{\sa} \ot\,\e$ to this equation, we get that
the element $\s(U^\a_{ij})$ is equal to $\sum_{k=1}^{N_\a} \k^{-2}(U^\a_{ik}) \, \e(\s(U^\a_{kj})) \in \sa_\a$. Hence,
$\s(\sa_\a) \subseteq \sa_\a$. Define the bijective linear operator $P_\a$ on $\sa_\a$ by setting $P_\a(a) = \s(a)$,
for all $a \in \sa_\a$. By assumption, $\ip{P_\a \, a}{a}
 = h(a^* \s(a)) \geq 0$, for all $a \in A_\a$; hence, $P_\a$
is positive. This allows us to define complex powers of
$P_\a$.

Now define the function $\gamma$ from $\C$ into the space of linear maps on $\sa$  such that for $z \in \C$ and $\a \in
\hat{G}$, we have that $\gamma_z(\sa_\a) \subseteq \sa_\a$ and $\gamma_z(a) = P_\a^{-iz}\,a$, for all $a \in \sa_\a$.
It is clear that  $\gamma_{ni} = \s_{\!\sa}^n$, for all $n \in \N_0$.

By diagonalising each operator $P_\a$, we can find a homogeneous orthonormal basis $(e_j)_{j \in J}$ for $\sa$ and
positive numbers $(\lambda_j)_{j \in J}$ such that $\s(e_j) = \lambda_j\, e_j$, for all  $j \in J$. Then $\gamma_z(e_j)
= \lambda_j^{-iz}\,  e_j$,  for all $j \in J$ and $z \in \C$. Hence, Theorem~\ref{one.prop1} implies that $\s_{\!\sa}$
is induced by an exponential one-parameter group, which must be equal to $\gamma$.

Let $z \in \C$. We will extend $\gamma_z$ to $\Om$. Let $m$
and  $n$ be a nonnegative and a positive integer,
respectively, and let $v_{jk}$ be elements in $\sa$, where
${j=1,\dots,m}$ and ${ k=0,\ldots,n}$, such that
${\sum_{j=1}^m v_{j0} \, dv_{j1} \ldots dv_{jn} = 0}$. Define
the map ${f : \C \rightarrow \Om}$ by setting $f(y) =
\sum_{j=1}^m \gamma_y(v_{j0}) \, d(\gamma_y(v_{j1})) \ldots
d(\gamma_y(v_{jn}))$, for all $y \in \C$. The function $f$ is
clearly of finite exponential type. Moreover, for $l \in \Z$,
we have
\begin{eqnarray*}
f(l i) & = & \sum_{j=1}^m \gamma_{li}(v_{j0}) \,
d(\gamma_{li}(v_{j1})) \ldots d(\gamma_{l i}(v_{jn}))
\\ & = & \sum_{j=1}^m \s^l(v_{j0}) \, d(\s^l(v_{j1}))
\ldots d(\s^l(v_{jn})) = \s^l\bigl(\,\sum_{j=1}^m v_{j0} \, dv_{j1} \ldots dv_{jn}\,\bigr) = 0 \ .
\end{eqnarray*}
Hence $f = 0$. In particular, ${f(z) = \sum_{j=1}^m
\gamma_z(v_{j0}) \, d(\gamma_z(v_{j1})) \ldots
d(\gamma_z(v_{jn})) = 0}$.

This allows us to define a map $\b$ from $\C$ into the space
of linear maps on $\Om$  such that each $\b_z$ is of degree
zero and
\begin{equation} \label{one.eq1}
\b_z(v_0 dv_1 \ldots dv_n) = \gamma_z(v_0) \, d(\gamma_z(v_1))
\ldots d(\gamma_z(v_n)),
\end{equation}
for all $z \in \C$, $n \in \N_0$ and $v_0,\ldots,v_n \in \sa$.
Hence, each $\b_z$ is an algebra homomorphism extending
$\gamma_z$ and  $\b_z \,d = d\, \b_z$. The equality
$\gamma_z(a)^* = \gamma_{\bar{z}}(a^*)$, for all $a \in \sa$,
extends easily to the equality $\b_z(\om)^* =
\b_{\bar{z}}(\om^*)$, for all $\om \in \Om$.

Furthermore, because $\gamma_y\,\gamma_z = \gamma_{y+z}$, we
get $\b_y\,\b_z = \b_{y+z}$, for all $y,z \in \C$. Therefore,
$\b$ is a one-parameter group on $\Om$. Combining
Equation~(\ref{one.eq1}) with the fact that $\gamma$ is an
exponential one-parameter group on~$\sa$
 shows that $\b$ is exponential. Since
$\int \s = \int$, the remarks preceding
Theorem~\ref{one.prop1} imply that $\int \b_z = \int$.

If $\om \in \Om$, the functions $z \mapsto (\tau_z \ot
\b_z)\D_\Om(\om)$ and ${z \mapsto \D_\Om(\b_z(\om))}$ from
$\C$ to ${\sa \ot \Om}$ are of finite exponential type. For $n
\in \N_0$, we have $(\tau_{ni} \ot \b_{ni})\D_\Om(\om) =
(\k^{-2n} \ot \s^n)\D_\Om(\om) = \D_\Om(\s^n(\om)) =
\D_\Om(\b_{ni}(\om))$. Therefore, the two functions agree on
the whole complex plane. We conclude that $(\tau_z \ot
\b_z)\D_\Om = \D_\Om \b_z$, for all $z \in \C$. Clearly, this
equality implies that $\b_z(\Om^{\text{inv}}) =
\Om^{\text{inv}}$, for all $z \in \C$.~\qed

The condition $h(a^* \s(a)) \geq 0$ is not easy to check. However,  in the next result we give some equivalent
conditions, of which the third one can be easily checked on an algebraic level.

\begin{cor} \label{one.cor1} The following conditions are equivalent:
\begin{enumerate}
\item $h(a^* \s(a)) \geq 0$, for all $a \in \sa$;
\item There exists an algebra automorphism $\phi$ on $\sa$  such that
$\phi^2 = \s_{\!\sa}$, $h \phi =  h$ and $\phi(\phi(a^*)^*) = a$, for all $a \in \sa$;
\item There exists an algebra automorphism $\phi$ on $\sa$ such that
$\phi^2 = \s_{\!\sa}$, \newline $(\tau_{i/2} \ot \phi)\D = \D \phi$ and $\phi(\phi(a^*)^*) = a$, for all $a \in \sa$
\end{enumerate}
\end{cor}

\demo{Proof} First we show $(1) \Rightarrow (3)$. Therefore, assume Condition~1. By the theorem, there exists a unique
exponential one-parameter group $\b$ on $\sa$ such that $\b_i = \s_{\!\sa}$. Put $\phi = \b_{i/2}$. Then $\phi$
satisfies the conditions in the Condition~3.

The implication $(3) \Rightarrow (2)$ follows from fact that
the equality $(\tau_{i/2} \ot \phi)\D = \D \phi$ implies $h
\phi = h$.

Finally, we show the implication $(2) \Rightarrow (1)$ holds: Assuming (2), we get, for $a \in \sa$, $ h(a^* \s(a))
=h(a^* \phi^2(a)) = h\bigl(\phi^{-1}(a^* \phi^2(a))\bigr) = h(\phi^{-1}(a^*) \phi(a)) = h(\phi(a)^* \phi(a)) \geq
0$.~\qed

It is useful to note that, of course, the equations in
Condition~3 of the corollary  hold for all $a$ in $\sa$ if
they hold for all $a$ in a set of algebra generators of~$\sa$.

Suppose that $h(a^* \s(a)) \geq 0$, for all $a \in \sa$. Let $\b$ denote the unique exponential one-parameter group on
$\Om$ such that $\b_i = \s$. We denote by $\g$ and $\g'$ the exponential one-parameter groups on $\sa$ obtained,
respectively, by restricting $\b$  and by setting $\gamma'_z = R \, \gamma_z\,R$, for all $z \in \C$.

Since $\s_\sa\,\rho_z = \rho_z \, \s_{\!\sa}$ and $\s_{\!\sa}\,\tau_z = \tau_z \, \s_{\!\sa}$, for all $z \in \C$,
Lemma~\ref{one.lem1} implies that $\g_y\,\rho_z = \rho_z \, \g_y$ and $\g_y \, \tau_z = \tau_z \, \g_y$, for all $y,z
\in \C$. Because $\k = R \tau_{-i/2} = \tau_{-i/2} R$, we get $ \gamma_z' = \k \gamma_z \, \k^{-1} = \k^{-1} \gamma_z
\, \k$. Also, $\gamma_y'\,\rho_z = \rho_z \, \gamma_y'$ and $\gamma_y' \, \tau_z = \tau_z \, \gamma_y'$, since
$R\r_z=\r_{-z}R$ and $R\t_z=\t_zR$.

\begin{lem} \label{lem: gamma-lemma} Suppose that
$h(a^*\s(a))\ge 0$, for all $a\in\sa$. Let $y,z \in \C$. Then
\begin{enumerate}
\item $(\gamma_z' \ot \tau_z)\D = \D \, \gamma_z'$;
\item $(\gamma_z' \ot \gamma_z)\D = \D \, \gamma_z' \,
\gamma_z\,\tau_{-z}$;
\item $\gamma_z' \, \gamma_y = \gamma_y \, \gamma_z'$.
\end{enumerate}
\end{lem}

\demo{Proof} Theorem~\ref{one.prop2} implies that $(\tau_z \ot
\gamma_z)\D = \D \gamma_z$. Hence, if $\flip$ is the flip map,
\begin{eqnarray*}
& & (\gamma_z' \ot \tau_z)\D = (\k \gamma_z \, \k^{-1} \ot \k
\tau_z \, \k^{-1})\D = (\k \gamma_z \ot \k \tau_z)\flip\D
\k^{-1}
\\ & & \spat = \flip (\k \tau_z \ot  \k \gamma_z)\D \k^{-1}
= \flip(\k \ot \k)\D \gamma_z \, \k^{-1} = \D \k
\gamma_z\,\k^{-1} = \D \gamma_z',
\end{eqnarray*}
and therefore,
\begin{eqnarray*}
& & (\gamma_z' \ot \gamma_z) \D = (\gamma_z'\,\tau_z \ot
\gamma_z \, \tau_z)\D \tau_{-z} = (\gamma_z' \, \tau_z \ot
\tau_z \, \gamma_z)\D \tau_{-z}
\\ & & \spat = (\gamma_z' \ot \tau_z)\D \gamma_z\, \tau_{-z} = \D
\gamma_z'\,\gamma_z\,\tau_{-z} \ .
\end{eqnarray*}
Also,
\begin{eqnarray*}
& & \D \gamma_z' \, \gamma_y =  (\gamma_z' \ot \tau_z)\D\g_y =
(\gamma_z' \, \tau_y \ot \tau_z \, \gamma_y)\D
\\ & & \spat = (\tau_y \, \gamma_z' \ot \gamma_y \, \tau_z)\D = (\tau_y \ot
\gamma_y)\D \gamma_z' = \D \gamma_y\,\gamma_z' \ .
\end{eqnarray*}
Therefore, injectivity of $\D$ implies that
$\gamma_z'\,\gamma_y = \gamma_y\,\gamma_z'$. \qed

\begin{lem} \label{one.prop3} Suppose that
$h(a^*\s(a))\ge 0$, for all $a\in\sa$. For every  integer $k \in \N_0$, there exists a basis $(\om_{k n})_{n=1}^{m_k}$
for $\ominv_k$ and positive numbers $(\nu_{kn})_{n=1}^{m_k}$ such that $\b_z(\om_{kn}) = \nu_{kn}^{-iz} \, \om_{kn}$,
for all $n = 1,\ldots,m_k$ and $z \in \C$.
\end{lem}

\demo{Proof} Define the linear map $w : \sa \rightarrow \Om_1$ by setting $w(a) = m(\k \ot d)\D(a)$, for all $a \in
\sa$. Here $m : \sa \ot \Om_1 \rightarrow \Om_1$ is the linear map given by $m(a \ot \om) = a \om$, for all $a \in \sa$
and $\om \in \Om_1$. Observe that $w(\sa) = \ominv_1$. We refer the reader to \cite[Section~14.1.1]{KS} for more
information on the map $w$.


By Lemma~\ref{lem: gamma-lemma}, we have, for all $a \in \sa$,
\begin{eqnarray}
 \hspace{-2ex} \b_z(w(a)) \hspace{-1ex}  & = & \hspace{-1ex} \b_z\bigl( m(\k \ot d)\D(a)\bigr) = m(\gamma_z \, \k \ot \b_z\, d)\D(a) \nonumber  = m(\k
\gamma_z' \ot d\, \gamma_z)\D(a) \\ &  = & \hspace{-1ex} m(\k
\ot d)\D((\gamma_z' \, \gamma_z \, \tau_{-z})(a)) =
w((\gamma_z' \, \gamma_z \, \tau_{-z})(a)) \ . \label{one.eq2}
\end{eqnarray}

Define the automorphism $\phi$ on $\sa$ by setting $\phi =
\gamma_{i/2}' \, \gamma_{i/2} \, \tau_{-i/2}$. Then $\phi^2 =
\gamma_i' \, \gamma_i \, \tau_{-i}$, $h \, \phi = h$ and
$\phi(\phi(a^*)^*) = a$, for all $a \in \sa$. As in the proof
of Corollary~\ref{one.cor1}, this implies that $h(a^*
\phi^2(a)) \geq 0$, for all $a \in \sa$.

Let $\a \in \hat{G}$. In the same way as in the beginning of the proof of Theorem~\ref{one.prop2}, the equalities $\D
\gamma_{i/2}' = (\gamma_{i/2}' \ot \tau_{i/2})\D$ and $\tau_{i/2}(\sa_\a) =\sa_\a$ allow us to conclude that
$\gamma_{i/2}'(\sa_\a) =\sa_\a$. Similarly, $\gamma_{i/2}(\sa_\a) =\sa_\a$. Hence, $\phi(\sa_\a)= \sa_\a$, for all $\a
\in \hat{G}$. We can therefore use the same idea as in the proof of Theorem~\ref{one.prop2} to find a basis $(f_j)_{j
\in J}$ of $\sa$ and positive numbers $(\mu_j)_{j \in J}$ such that $\phi^2(f_j) = \mu_j\,f_j$, for all $j \in J$.

Equation~(\ref{one.eq2}) gives $\b_i(w(f_j)) = \mu_j \,
w(f_j)$, for all $j\in J$. Hence, $\b_{li}(w(f_j)) = (\b_i)^l
(w(f_j)) = \mu_j^l \, w(f_j)$, for all $l \in \N_0$. The
functions ${z \mapsto \b_z(w(f_j))}$ and ${z \mapsto
\mu_j^{-iz} \, w(f_j)}$ from $\C$ to $\Om_1$ are both of
finite exponential type and agree on $\N_0\, i$, implying that
they are equal on the whole complex plane. Thus, $\b_z(w(f_j))
= \mu_j^{-iz} \, w(f_j)$, for all $z \in \C$.

The lemma now follows easily, using the fact that
$\Om^{\text{inv}}_k$ is the linear span of the elements
$w(f_{j_1}) \ldots w(f_{j_k})$, where $j_1,\ldots,j_k \in J$,
combined with the fact that, for $z \in \C$, the
multiplicativity of $\b_z$ implies that $\b_z\bigl(w(f_{j_1})
\ldots w(f_{j_k})\bigr) = (\mu_{j_1} \ldots \mu_{j_k})^{-iz}$
$w(f_{j_1}) \ldots w(f_{j_k})$, for all $j_1,\ldots,j_k \in
J$. \qed

By defining the inner-product on each space $\ominv_k$ in such
a way that the basis in the previous lemma is orthonormal, we
get the following theorem.

\begin{thm} \label{thm: inner product theorem}
There exists a graded, left-invariant inner product $\ip{\cdot}{\cdot}$ on $\Om$ such that $\ip{\s(\om)}{\om} \geq 0$,
for all $\om \in \Om$, if, and only if, $h(a^*\s(a))\ge 0$, for all $a\in\sa$.
\end{thm}

\begin{thm} \label{thm: inner product theorem 2} Suppose given a graded, left-invariant inner
product $\ip\cdot\cdot$ on~$\Om$ relative to which $\s$ is
positive. Let $\b$ be the unique exponential one-parameter
group on $\Om$ such that $\b_i=\s$. Then the following
conditions hold:
\begin{enumerate}
\item $\ip{\b_z(\om)}{\om}  \geq 0$, for all $\om \in \Om$ and
$z \in \R\,i$.
\item $\ip{\b_z(\om)}{\eta} = \ip{\om}{\b_{-\bar{z}}(\eta)}$,
 for all $\om, \eta \in \Om$.
\item $d^* \, \b_z = \b_z \, d^*$, $D \, \b_z = \b_z \, D$ and
$\nab \, \b_z = \b_z \, \nab$, for all $z \in \C$.
\item For $L$ the Hodge operator, we have $L \, \b_z = \b_z \,
L$,
for all $z \in \C$.
\end{enumerate}
\end{thm}
\demo{Proof} Combining Theorem~\ref{one.prop2}
and~Lemma~\ref{one.prop3} shows that there exists an
orthonormal basis $(\eta_k)_{k \in K}$ of $\Om$ and a family
of positive numbers $(\nu_k)_{k \in K}$ such that
$\b_z(\eta_k) = \nu_k^{-iz}\,\eta_k$, for all $k \in K$. From
this, Conditions~1 and~2 follow. Condition~2 implies that, if
$z \in \C$, then $\b_{-\bar{z}}$ is adjointable and that its
adjoint is $\b_z$. Hence, taking the adjoint of the equation
$d \,\b_{-\bar{z}} = \b_{-\bar{z}} \, d$, we get $\b_z\,d^* =
d^* \,\b_z$. The other two equalities of Condition~3 follow
immediately.

If $\om,\eta \in \Om$, then, using the fact that $\int \s =
\int$, we have
\begin{eqnarray*}
\int \om^* L(\s(\eta))  & = & \ip{\s(\eta)}{\om} =
\ip{\eta}{\s(\om)} = \int \s(\om)^* L(\eta)
\\ & = & \int \s(\s(\om)^*)\,\s(L(\eta)) = \int \om^*\, \s(L(\eta)) \ .
\end{eqnarray*}
Hence, by  faithfulness of $\int$, $L(\s(\eta)) = \s(L(\eta))$
and therefore, by Lemma~\ref{one.lem1} (more precisely, by a
result for $\Om$ that is analogous), $L\, \b_z = \b_z \,L$,
for all $z \in \C$. Thus, Condition~4 holds.\qed

\begin{thm} Suppose that $h(a^*\s(a))\ge 0$, for all
$a\in\sa$. Then there exists a graded, left-invariant inner
product $\ip{\cdot}{\cdot}$ on~$\Om$ such that
$\ip{\s(\om)}{\om} \geq 0$, for all $\om \in\Om$, and such
that the Hodge operator $L$ associated to this inner product
satisfies the conditions of Theorem~\ref{thm: L-properties
theorem}.\end{thm}

\demo{Proof} First we choose, as we may, a graded,
left-invariant inner product $\ip{\cdot}{\cdot}'$ on $\Om$
such that $\ip{\s(\om)}{\om}' \geq 0$, for all $\om \in \Om$.
Now construct a new graded, left-invariant inner product
$\ip\cdot\cdot$ by the method outlined in Section~\ref{sec:
inner product}. One easily checks that $\ip{\s(\om)}\om\ge 0$,
for all $\om\in \Om_k$, if $k\ne N/2$, by using the fact that
$\s$ and the Hodge operator $L'$ associated to
$\ip\cdot\cdot'$ commute. We show now that by refining the
method used in Section~\ref{sec: inner product} to define the
inner product on $\Om_k$, when $k=N/2$, we can arrange to have
in this case also that $\ip{\s(\om)}\om\ge 0$, for all $\om\in
\Om_k$.

Let $T$ be the unique linear operator on~$\Om_k^{\text{inv}}$
such that $\ip{T \om}{\eta}' = \int \eta^* \om$, for all
$\om,\eta \in \Om^{\text{inv}}_k$. Then $T^* =
(-1)^{k(N-k)}\,T$. Using the fact that $\s$ is a positive
operator with respect to $\ip{\cdot}{\cdot}'$, and that, for
$\om,\eta \in \Om_k^{\text{inv}}$, $\int \eta^* \s(\om) = \int
\s^{-1}(\eta^* \s(\om)) = \int \s(\eta)^* \, \om$, it follows
that the operators $\s$ and $T$ commute. Since these operators
are normal, there is an orthogonal basis $u_1,\ldots,u_m$ of
$\Om^{\text{inv}}_k$ consisting of eigenvectors of $\s$ and
$T$. Setting $\e_k=1$, if $(-1)^{k(N-k)}=1$ and $\e_k=i$, if
$(-1)^{k(N-k)}=-1$, we can re-normalize the basis elements
$u_j$ to get a new basis ${e_1,\dots,e_m}$ of eigenvectors of
$\s$ such that ${\int e^*_qe_p}=\pm\e_k\d_{pq}$, for
${p,q=1,\dots,m}$. We now proceed as before, in
Section~\ref{sec: inner product}, to choose a new inner
product on $\ominv_k$ making the basis $e_j$ orthonormal. From
this we complete our construction to get a graded,
left-invariant inner product $\ip\cdot\cdot$ on $\Om$. For
this inner product the associated Hodge operator $L$ satisfies
the conditions of Theorem~\ref{thm: L-properties theorem}.
Also, $\s$ is clearly positive for this inner product.~\qed

\section{Woronowicz's 3-dimensional calculus}\label{sec: eigenvalues
of nabla}

In this section we compute the eigenvalues of the Laplacian
for a special 3-dimensional calculus over the quantum group
$SU_q(2)$, where $q$ is a real number for which $0<\m q < 1$.
This calculus has a special interest, since it was the first
calculus constructed by Woronowicz in the quantum group
setting.

We begin by recalling some basic facts about the calculus
(see~\cite{WOR1} for details). It is a left-covariant,
3-dimensional $*$-differential calculus $\omd$ over the Hopf
algebra~$\sa_q$ underlying $SU_q(2)$. The space $\ominv_1$ has
a linear basis $\eta_1,\eta_2, \eta_3$ such that
$\eta_1\eta_2,\eta_2\eta_3,\eta_3\eta_1$ is a basis of
$\ominv_2$ and $\eta_1\eta_2\eta_3$ is one of $\ominv_3$. The
involution on $\Om$ is determined by the relations
\[
\eta^*_1=q\eta_3,\ \eta^*_2=-\eta_2,\ \eta^*_3=q^{-1}\eta_1.\]

Left-covariance implies that there exists linear functionals
$\chi_1,\chi_2,\chi_3$ on $\sa$ such that
\[ d a=\sum_{r=1}^3 (\chi_r * a) \eta_r,\]
for all $a\in \sa$. We note that \begin{equation} \label{eqn:
chi* equations}\chi_1^*=-q^{-1}\chi_3,\ \chi_2^*=\chi_2,\
\chi^*_3=-q\chi_1.\end{equation} As usual, we shall set
$\chi={\sum_{r=1}^3 \chi^*_r\chi_r}$. Note that
$d(\ominv_2)=0$, by Equation Set~(3.13) of~\cite{WOR1}.

The following commutation relations will be needed frequently
\cite[Table 5]{WOR1} :
\[ \eta_2\eta_1=-q^4\eta_1\eta_2,\
\eta_3\eta_2=-q^4\eta_2\eta_3,\
\eta_1\eta_3=-q^{-2}\eta_3\eta_1.\]

Note also that $\eta^2_r=0$, for $r=1,2,3$. It follows that
$\eta_r\eta_s\eta_t=c\eta_1\eta_2\eta_3$, for some scalar $c$
and that $c\ne 0$ if, and only if, the subscripts $r,s,t$ are
distinct.

We set $\ta=\eta_1\eta_2\eta_3$. Then $\ta^*=\ta$.

We now choose an inner product on $\ominv$ such that the
vectors
\[1\,,\ \eta_1\,,\ \eta_2\,,\ \eta_3\,,\ \eta_1\eta_2\,,\ \eta_2\eta_3\,,\ \eta_3\eta_1\,,
\ \eta_1\eta_2\eta_3\] form an orthonormal basis. In the usual way, this allows us to define a graded, left-invariant
inner product on~$\Om$.  Then $\ta$ is a volume element, since $\ta^*=\ta$ and $\n \ta=1$. We let $\int$ denote the
integral associated to $\ta$.

If $\eta_r\eta_s\eta_t=c\eta_1\eta_2\eta_3$, then ${c=\int
\eta_r\eta_s\eta_t}$; hence, ${\int \eta_r\eta_s\eta_t\ne0}$
if, and only if, $r,s,t$ are distinct.

\begin{thm} $\Om$ is non-degenerate. \end{thm}

\demo{Proof} We have to show that if $0\le k\le 3$ and $\om\in
\Om_k$ and $\om'\om=0$, for all $\om'\in \Om_{3-k}$, then
$\om=0$.

Consider first the case where $k=0$, so that $\om=a$, for some element $a\in \sa_q$. Then, we suppose $\ta b a=0$, for
all $b\in\sa_q$. If we take $b=1$, this becomes $\ta a=0$ and therefore, $a=0$, as required.

Now we consider the case where $k=1$. In this case we may write $\om=\sum_{r=1}^3 \eta_r a_r$, for some elements
$a_1,a_2,a_3\in \sa_q$. If $s,t$ are distinct elements of $\{1,2,3\}$, then $0=\eta_s\eta_t\om=\eta_s\eta_t\eta_r a_r$,
where $r,s,t$ are distinct. Hence, $a_r=0$. It follows that $\om=0$.

Suppose now $k=2$. In this case, we may write
$\om=\eta_1\eta_2a_1+\eta_2\eta_3a_2+\eta_3\eta_1a_3$, for
some elements $a_1,a_2,a_3\in \sa_q$. Then
$0=\eta_3\om=\eta_3\eta_1\eta_2a_1$. Hence, $a_1=0$.
Similarly, $a_2=0$ and $a_3=0$ and therefore, $\om=0$.

If $k=3$, then $\om=1\om=0$.

Thus, in all cases, $\om=0$ and therefore $\Om$ is
non-degenerate, as required. \qed

It follows from non-degeneracy  of $\Om$ that $\int$ is
faithful (Lemma~\ref{thm: integral}) and a Hodge operator $L'$
exists (Theorem~\ref{thm: existence of Hodge operator}).  We
need to know the values of $L'$ on $\eta_1,\eta_2,\eta_3$.
First, set $L'(\eta_1)=
\l_1\eta_1\eta_2+\l_2\eta_2\eta_3+\l_3\eta_3\eta_1$, where
$\l_1,\l_2,\l_3$ are scalars to be determined. Observe that
${\int\eta^*_sL'(\eta_r)}=\ip{\eta_r}{\eta_s}=\d_{rs}$. Hence,
${1=\int \eta^*_1L'(\eta_1)}={\int
q\eta_3\l_1\eta_1\eta_2}={q\l_1\int (-q^2)\eta_1\eta_3\eta_2}=
{-q^3\l_1\int \eta_1(-q^4\eta_2\eta_3)}={q^7\l_1\int
\eta_1\eta_2\eta_3=q^7\l_1}$. Hence, $\l_1=q^{-7}$. Similar
reasoning shows that $\l_2=0$ and $\l_3=0$. Continuing this
way we can calculate that
\[ L'(\eta_1)=q^{-7}\eta_1\eta_2,\
L'(\eta_2)=-q^{-6}\eta_3\eta_1,\ L'(\eta_3)=q\eta_2\eta_3.\]

Now we redefine the inner product on $\Om$ as in
Section~\ref{sec: inner product} to obtain a new Hodge
operator $L$ satisfying the conditions of Theorem~\ref{thm:
L-properties theorem}. The new inner product remains unchanged
on $\ominv_0$, $\ominv_1$ and $\ominv_3$ (since $L'(1)=\ta$),
but on $\ominv_2$ it is chosen in such a way that
$L$ maps $\ominv_1$ isometrically onto
$\ominv_2$. Hence, for the new inner product, the elements
\[1,\eta_1,\eta_2,\eta_3,q^{-7}\eta_1\eta_2,q\eta_2\eta_3,
-q^{-6}\eta_3\eta_1, \eta_1\eta_2\eta_3\] form
an orthonormal basis.

Since $L=L'$ on $\ominv_1$ and since $L^2=1$ we get
\[ L(\eta_1)=q^{-7}\eta_1\eta_2,\
L(\eta_2)=-q^{-6}\eta_3\eta_1,\ L(\eta_3)=q\eta_2\eta_3 \] and
\[ L(\eta_1\eta_2)=q^7\eta_1,\ L(\eta_2\eta_3)=q^{-1}\eta_3,\
L(\eta_3\eta_1)=-q^6\eta_2.\]

Now set $E_0=\id_{\sa_q}$ and $E_r=E_{\chi_r}$, for $r=1,2,3$.
Also, denote by $M_0$ the restriction of $d$ to $\ominv$ and,
for $r=1,2,3$, denote by $M_r$ the operator on $\ominv$ of
left multiplication by $\eta_r$. Using the usual
identification of $\Om$ with ${\sa_q\otimes \ominv}$, we have
${d=\sum_{r=0}^3 E_r\otimes M_r}$, by Equation~(\ref{eqn:
tensor d-equation}). If $\ell$ denotes the restriction of $L$
to $\ominv$, then clearly ${L=\id_{\sa_q}\otimes \ell}$. Hence
$L d= {\sum_{r=0}^3 E_r\otimes \ell M_r}$.

Now let $T$ be the restriction map ${Ld\colon \Om_1\to
\Om_1}$, and, for $r=0,1,2,3$, let $T_r$ be the restriction of
$\ell M_r$ to $\ominv_1$. Then ${T=\sum_{r=0}^3 E_r\otimes
T_r}$.

Observe that $T=T^*$ and that $T^2(\om)=d^*d\om$, for all
$\om\in \Om_1$. For, if ${\om,\om'\in \Om_1}$, then
$\ip\om{T\om'}=\ip\om{Ld\om'}=\ip\om{d^*L\om'}$, by
Theorem~\ref{thm: L-properties theorem}. Hence,
$\ip\om{T\om'}=\ip{d\om}{L\om'}=\ip{Ld\om}{L^2\om'}=\ip{T\om}{\om'}$,
since $L^2=\id$. Hence, $T$ is adjointable, with adjoint
$T^*=T$. Also, $T^2(\om)=L d L d(\om)=d^*L L d(\om)
=d^*d(\om)$.

This identification of a self-adjoint square root $T$ for the
restriction of $d^*d$ to $\Om_1$ will be the key to our
approach to calculating the eigenvalues of the Laplacian
$\nabla$ on $\Om_1$. In fact, we shall calculate the
eigenvalues of $T$ and use these and some calculations from
Section~\ref{sec: 0-forms diagonalisation} to find the
eigenvalues of $\nabla$ on~$\Om_1$.

It is straightforward to calculate the matrices of the
operators $T_0,T_1,T_2,T_3$ (relative to the basis
$\eta_1,\eta_2,\eta_3$), using the following Cartan--Maurer
formulas from~\cite[Table 6]{WOR1}:
\[ d\eta_1=q^2(1+q^2)\eta_1\eta_2,\ d\eta_2=q\eta_1\eta_3,\
d\eta_3=q^2(1+q^2)\eta_2\eta_3.\] Hence, $
Ld\eta_1=q^9(1+q^2)\eta_1$, $Ld\eta_2=q^5\eta_2$ and
$Ld\eta_3= q(1+q^2)\eta_3$. Therefore, the matrix form of
$T_0$ is given by
\[T_0=\left[\begin{array}{ccc} q^9(1+q^2)&0&0\\
                     0&q^5&0\\
                     0&0&q(1+q^2)\end{array}\right].\]
Since $T_1(\om)=L(\eta_1\om)$, we have $T_1(\eta_1)=0$,
$T_1(\eta_2)=q^7\eta_1$ and  $T_1(\eta_3)=q^4\eta_2$, which
gives us the matrix of~$T_1$. We can calculate $T_2$ and $T_3$
similarly. We get
\[T_1=\left[\begin{array}{ccc} 0&q^7&0\\
                     0&0&q^4\\
                     0&0&0\end{array}\right],\
T_2=\left[\begin{array}{ccc} -q^{11}&0&0\\
                     0&0&0\\
                     0&0&q^{-1}\end{array}\right],
\ T_3=\left[\begin{array}{ccc} 0&0&0\\
                     -q^6&0&0\\
                     0&-q^3&0\end{array}\right].\]
Now, if $\eta_{rs}$ is the system of matrix units for
$B(\ominv_1)$ corresponding to the basis
$\eta_1,\eta_2,\eta_3$, then $T={\sum_{r,s=1}^3 T_{rs}\otimes
\eta_{rs}}$, where, from our calculations,
\[ \left[\begin{array}{ccc} T_{11}&T_{12}&T_{13}\\
                     T_{21}&T_{22}&T_{23}\\
                     T_{31}&T_{32}&T_{33}\end{array}\right]=
 \left[\begin{array}{ccc} q^9(1+q^2)-q^{11}E_2&q^7E_1&0\\
                     -q^6E_3&q^5&q^4E_1\\
                     0&-q^3E_3&q(1+q^2)+q^{-1}E_2\end{array}\right].\]
It is useful to observe that, if $a\in \sa_q$, then
\[T(a\eta_s)=\sum_{r=1}^3T_{rs}(a)\eta_r.\]

The formula for $T$ we have derived is still not explicit
enough to enable us to calculate eigenvalues. To make things
fully explicit, we now need to use some aspects of the
representation theory of $SU_q(2)$. Recall that this quantum
group has a complete system of inequivalent invertible
irreducible representations ${W^1,W^2,\dots, W^M,\dots}$,
where $W^M$ acts on a Hilbert space $H_M$ of dimension $M$.

The $W^M$ are defined as in \cite[Equation~5.32]{WOR1}. They
are not unitary, which is not in line with our practice up to
this of only dealing with unitary representations. However, we
shall be making extensive use of various formulas from
\cite{WOR1}, so we use the $W^M$. There is a standard way of
constructing an equivalent unitary representation to a given
invertible representation. If $Q_M=h(W^{M*}W^M)$ (that is, the
scalar matrix got by applying $h$ to each entry of
$W^{M*}W^M$), then $Q_M$ is positive and invertible
and $U^M=Q_M^{1/2}W^MQ_M^{-1/2}$ is an irreducible
unitary representation equivalent to $W^M$. See
\cite[Theorem~5.2]{WOR2} for details.

We shall need some quasi-orthogonality relations for the
$W^M$. Since the matrix entries of $U^M$ are orthogonal to
those of $U^N$, if $M\ne N$, the matrix entries of $W^M$ are
orthogonal to those of $W^N$.

It is well known that $U^M_{pk}$ and $U^M_{tu}$ are orthogonal if $k\ne u$ \cite[Theorem~7.3]{MT}.  It turns out that
the elements $W^M_{pk}$ have this property also and this will be important in the sequel. The reason is that in this
special case the matrix $Q_M$ is diagonal, so that there are positive numbers $q_M(k)$ such that
$U^M_{pk}=q_M(p)q_M(k)^{-1}W^M_{pk}$.


For this example we even have that $h((W^M_{pk})^* W^M_{p'l}) = 0$ if $k \not=l$ or $p \not=  p'$. A proof of this fact
will be given in Lemma \ref{lem:orthogonality}.

\medskip

Fix $M\ge 1$. Set $A_r=\chi_r(W^M)$, for $r=1,2,3$ (these
matrices are the {\em infinitesimal generators} of the
representation~$W^M$, see \cite[Theorem~5.4]{WOR1}). Let $m$
be the integer or half-integer for which $M=2m+1$. Then the
space $H_M$ admits an orthogonal basis ${\xi_{-m},\xi_{-m+1},
\dots, \xi_{m-1},\xi_m}$ for which we have
\begin{equation} \label{eqn: A-equation}
A_1\xi_k=-c_{k+1}\xi_{k+1},\ A_2 \xi_k=\l_k\xi_k,\
A_3\xi_k=qc_k\xi_{k-1}.\end{equation}
Here, for arbitrary $k$,
we set
\[ \l_k=q^2(1-q^2)^{-1}(q^{-4k}-1)\]
and, for ${-m\le k\le +m}$, we set
\[c_k=q(1-q^2)^{-1}[(q^{-2k}-q^{2m})(q^{-2m}-q^{-2(k-1)})]^{1/2}.\]
(Of course, $c_k$ is dependent on $m$, that is, on $M$, as
well as on~$k$, but we suppress this dependence in the
notation for the sake of simplifying it.)

So that our formulas will always make sense, we define
$\xi_{-m-1}=\xi_{m+1}=0$ and $c_{-m-1}=c_{m+1}=0$.

Denote by $W^M_{pk}$ the matrix entries of $W^M$ relative to
the basis ${\xi_{-m}, \dots,\xi_m}$. To calculate
$d(W^M_{pk})$ we need to calculate
$E_r(W^M_{pk})$, for $r=1,2,3$. We have \[E_1(W^M_{pk})=
{(\id\otimes \chi_1)\sum_{i=-m}^{i=+m} W^M_{pi}\otimes
W^M_{ik}}={\sum_{i=-m}^{i=+m}A_1(i
k)W^M_{pi}}={-c_{k+1}W^M_{p,k+1}}.\] By similar calculations
we get the following:
\begin{equation} \label{eqn: E-equations} E_1(W^M_{pk})=
-c_{k+1}W^M_{p,k+1},\ E_2(W^M_{pk})=\l_kW^M_{pk},\
E_3(W^M_{pk})=qc_kW^M_{p,k-1}.\end{equation} So that these
formulas always make sense, we define
$W^M_{p,-m-1}={W^M_{p,m+1}=0}$. For other purposes below, we
also set $W^M_{p,-m-2}=W^M_{p,m+2}=0$. It follows from the
preceding equations that for $k$ in the range ${-m,\dots,+m}$,
\begin{equation}\label{eqn: d-equation}
dW^M_{pk}=-c_{k+1}W^M_{p,k+1}\eta_1+\l_k W^M_{pk}\eta_2+ qc_kW^M_{p,k-1}\eta_3.\end{equation} Of course,
$TdW^M_{pk}=Ld^2W^M_{pk}=0$. Moreover, if $M>1$, then $dW^M_{pk}$ is non-zero. Otherwise, $W^M_{pk}=c1$, for some
scalar~$c$ (this is a ``connectedness'' result for $SU_q(2)$, see~\cite[Theorem~2.3]{WOR1}). Since $W^M_{pk}$ and the
unit 1 are orthogonal in $\sa_q$, because they are matrix elements of inequivalent representations of $SU_q(2)$, we
must have $c=0$. But then ${W^M_{pk}=0}$, an impossibility, since $W^M_{pk}$ belongs to a linear basis of $\sa_q$.
Hence, $dW^M_{pk}\ne 0$, as claimed. It follows that $dW^M_{pk}$ is an eigenvector for $T$.

From these observations, it is plausible to conjecture that it may be useful to look at the linear space spanned by the
vectors occurring in the righthand side of Equation~(\ref{eqn: d-equation}). Indeed, this turns out to be the key to
solving the eigenvalue problem for $T$.

For ${p=-m,-m+1,\dots,+m}$ and ${k=-m-1,-m,\dots,m+1}$, denote by $\EEONE$ the linear span of the vectors
$W^M_{p,k+1}\eta_1$, $W^M_{pk}\eta_2$ and $W^M_{p,k-1}\eta_3$. Since $\eta_1,\eta_2,\eta_3$ are orthogonal, so are
these spanning vectors. It follows that we have $\dim({\EEONE})=3$, if $k$ is in the range ${-m+1,\dots,m-1}$. If $M\ne
1$ and $k=-m$ or ${k=+m}$, $\dim(\EEONE)=2$. If $M\ne 1$ and $k=-m-1$ or ${k=m+1}$, then ${\dim(\EEONE)=1}$. Finally,
the spaces $\E_1(1,0,-1)$, $\E_1(1,0,0)$ and $\E_1(1,0,+1)$ all have dimension one also.

Set $f_1=W^M_{p,k+1}\eta_1$, $f_2=W^M_{pk}\eta_2$ and
$f_3=W^M_{p,k-1}\eta_3$.

We shall need the following result.

\begin{lem} \label{lem: d=0} Let $\om$ be an element of $\EEONE$ that is
orthogonal to $dW^M_{pk}$. Then $d^*\om=0$. \end{lem}

\demo{Proof}  We have $\om={\sum_{r=1}^3\a_r f_r}$, for some scalars $\a_1,\a_2,\a_3$. Hence, $d^*\om={\sum_{r=1}^3\a_r
d^*f_r}$. Using Equation Set~(\ref{eqn: E-equations}), and the relations $E^*_1=-q^{-1}E_3$, $E^*_2=E_2$ and
$E^*_3=-qE_1$, which follow from Lemma~\ref{thm: E-adjointability theorem} and Equation Set~(\ref{eqn: chi*
equations}), we calculate that $d^*f_1=-c_{k+1}W^M_{pk}$, $d^*f_2=\l_kW^M_{pk}$ and $d^*f_3=qc_kW^M_{pk}$. Hence,
$d^*\om={(-\a_1c_{k+1}+\a_2\l_k+\a_3qc_k)W^M_{pk}}$. However, by hypothesis
$\ip{d^*\om}{W^M_{pk}}={\ip{\om}{dW^M_{pk}}=0}$. Hence, $d^*\om=0$.~\qed

It is clear that the spaces $\EEONE$ are pairwise orthogonal; more precisely, $\EEONE \perp \E(M',p',k')$ if
$(M,p,k)\ne (M',p',k')$, since the elements $W^M_{pk}$ are orthogonal in $\sa_q$ and $\eta_1,\eta_2,\eta_3$ are
orthogonal in $\ominv_1$. Also, since $\Om_1$ is linearly spanned by the vectors $W^M_{pk}\eta_r$, it is the orthogonal
sum of the spaces $\EEONE$; that is,
\[\Om_1=\bigoplus_{(M,p,k)}\EEONE,\]
where $M$ ranges over the positive integers and
${p=-m,-m+1,\dots,+m}$ and ${k=-m-1,-m,\dots,m+1}$ (and
$M=2m+1$). We shall show that all the spaces occurring in this
decomposition of $\Om_1$ are invariant under~$T$. Hence, we
shall have reduced the eigenvalue problem for $T$ to the
problem of explicitly calculating $T$ on these spaces of
dimension~1, 2 or~3.

The relevance of this to the problem of finding the
eigenvalues of $\nabla$ can now be stated very explicitly. We
can split $\EEONE$ into the sum of $\mathbf{C}dW^M_{pk}$ and
its orthogonal complement ${\FF}$. By the same reasoning as we
used in Section~\ref{sec: 0-forms diagonalisation}, it follows
from Equation Set~(\ref{eqn: E-equations}) that $W^M_{pk}$ and
$dW^M_{pk}$ are eigenvectors of $\nabla$, with eigenvalue
$\nu_k$ given by \begin{equation}\label{eqn: nu-eigenvalue
equation} \nu_k= {c_{k+1}^2+\l^2_k+q^2c_k^2}.\end{equation}
(For $dW^M_{pk}$ to be an eigenvector, we need $M>1$.) Observe
that $\nu_k>0$, unless $m=k=0$, as is trivially verified. For
$\om\in \FF$, we have $d^*\om=0$, by Lemma~\ref{lem: d=0}.
Hence, $\nabla\om=d^*d\om=T^2\om$. Therefore, using the
invariance of the space $\EEONE$ under $T$ that we shall show
below, and that $T=T^*$, it follows from the equation
${TdW^M_{pk}=0}$ that $\FF$ is invariant under $T$ also and
therefore that $\FF$ is invariant under~$\nabla$; moreover,
the eigenvalues of $\nabla$ on this space are the squares of
the eigenvalues of $T$ on it. Thus, we have reduced the
eigenvalue problem for $\nabla$ to the problem of calculating
the eigenvalues of $T$ on the spaces $\FF$.

We calculate now the values of $T$ on the vectors
$W^M_{pk}\eta_1$, $W^M_{pk}\eta_2$ and $W^M_{pk}\eta_3$, where
$k$ is in the range ${-m,\dots,+m}$. We have
\begin{eqnarray*}
T(W^M_{pk}\eta_1)&=&{\sum_{r=1}^3 T_{r1}(W^M_{pk})\eta_r}\\&=&
{q^9(1+q^2)W^M_{pk}\eta_1-q^{11}E_2(W^M_{pk})\eta_1-
q^6E_3(W^M_{pk})\eta_2}
\\ &=&{(q^9(1+q^2)-q^{11}\l_k)W^M_{pk}\eta_1}
-{q^7c_kW^M_{p,k-1}\eta_2},\end{eqnarray*} by Equation
Set~(\ref{eqn: E-equations}). Similarly, one can show that
\[T(W^M_{pk}\eta_2)=-q^7c_{k+1}W^M_{p,k+1}\eta_1+q^{5}W^M_{pk}\eta_2
-q^4c_kW^M_{p,k-1}\eta_3\] and
\[T(W^M_{pk}\eta_3)=-q^4c_{k+1}W^M_{p,k+1}\eta_2+
q(1+q^2)W^M_{pk}\eta_3+q^{-1}\l_kW^M_{pk}\eta_3.\] However, it
is elementary to show that ${q^9(1+q^2)-q^{11}\l_k}=
-q^7\l_{k-1}$ and therefore also
${q(1+q^2)+q^{-1}\l_{k-1}}=q^3\l_k$. Hence,
\[T(f_1)=-q^7\l_kf_1-q^7c_{k+1}f_2,\]
\[T(f_2)=-q^7c_{k+1}f_1+q^5f_2-q^4c_kf_3\]
and
\[T(f_3)=-q^4c_kf_2+q^3\l_kf_3.\]
Hence, $\EEONE$ is invariant for $T$.

Consider now the case that $k$ is in the range
${-m+1,\dots,m-1}$, so that $f_1,f_2,f_3$ form an orthogonal
basis for $\E=\EEONE$. Our preceding calculations show
that the matrix form of the restriction $T_\E$ of $T$ relative
to this basis is given by
\[T_\E=\left[\begin{array}{ccc}
-q^7\l_k&-q^7c_{k+1}&0\\ -q^7c_{k+1}&q^5&-q^4c_k\\
0&-q^4c_k&q^3\l_k\end{array}\right].\] By Equation~(\ref{eqn:
d-equation}), $dW^M_{pk}$ belongs to $\E$ and we showed above
that $TdW^M_{pk}=0$. It therefore remains only to calculate
the eigenvalues of $T_\E$ on~$\F=\FF$.

Let us write
\[\left[\begin{array}{ccc}x&s&0\\
                          s&y&t\\
                          0&t&z\end{array}\right]=
T_\E=\left[\begin{array}{ccc} -q^7\l_k&-q^7c_{k+1}&0\\
-q^7c_{k+1}&q^5&-q^4c_k\\
0&-q^4c_k&q^3\l_k\end{array}\right].\] Then the characteristic
polynomial ${P(\l)=\det(T_\E-\l 1)}$ of $T_\E$ is given by
\[P(\l)={-\l^3+(x+y+z)\l^2-(xy+yz+xz-s^2-t^2)\l+xyz-t^2x-s^2z}.\]
Since 0 is an eigenvalue of $T_\E$, we must have ${x y
z-t^2x-s^2z=0}$. Hence,
\[P(\l)=-\l(\l^2-(x+y+z)\l+xy+yz+xz-s^2-t^2).\]
Set $B_k=x+y+z$ and $C_k=xy+yz+xz-s^2-t^2$. Then the
eigenvalues of $T$ on $\F$ are given by
\begin{equation}\label{eqn: mu-eigenvalue equation}
\mu_k^\pm=(B_k\pm\sqrt{B_k^2-4C_k})/2.\end{equation}
Explicitly, \begin{equation}\label{eqn: B-equation}
B_k=(q^3-q^7)\l_k+q^5\end{equation} and \begin{equation}
\label{eqn: C-equation}
C_k=q^8[(1-q^4)\l_k-q^2\l_k^2-q^6c_{k+1}^2-c_k^2]=
-q^{10}\nu_k.\end{equation} The second equality in
Equation~(\ref{eqn: C-equation}) follows from
Equation~(\ref{eqn: nu-eigenvalue equation}) and the first
equality in~(\ref{eqn: C-equation}) together with the easily
checked fact that ${\l_k=c^2_k-q^2c^2_{k+1}}$.

As observed above, the eigenvalues of $\nabla$ on $\FF$ are
the squares of those of $T$, that is, $(\mu^+_k)^2$ and
$(\mu^-_k)^2$.

We turn now to the case where $\E=\EEONE$ is 2-dimensional;
that is, where $M>1$ and $k=-m$  or $k=+m$.

Suppose first that $k=-m$. Then $f_3=0$ and $f_1,f_2$ form an
orthogonal basis of $\E$. Hence,
\[T_\E=\left[\begin{array}{cc}
-q^7\l_k&-q^7c_{k+1}\\ -q^7c_{k+1}&q^5\end{array}\right].\]
Since one eigenvalue of this matrix is zero (because
$TdW^M_{pk}=0$), the other eigenvalue on~$\F(M,p,-m)$ is the
sum of the diagonal entries ${q^5-q^7\l_k}$.

Similar considerations show that $T$ has  eigenvalue
$q^5+q^3\l_k$ on $\F(M,p,+m)$.

Finally, we consider the trivial case where $\EEONE$ is 1-dimensional: If $M>1$, then $k=-m-1$ or $k=m+1$, and it is
trivially checked that the eigenvalue of~$T$ on~$\EEONE$ is $-q^7\l_k$ or $ q^3\l_k$, respectively.  If $M=1$, then
$\EEONE=\E_1(1,0,k)=\F(1,0,k)$. Clearly, ${\F(1,0,-1)=\bC \eta_1}$ and we computed earlier that
${T\eta_1}={q^9(1+q^2)\eta_1}$; similarly, ${\F(1,0,0)=\bC \eta_2}$ and ${T(\eta_2)}={q^5\eta_2}$; finally, we have
${\F(1,0,+1)}={\bC\eta_3}$ and ${T\eta_3=q(1+q^2)\eta_3}$.

This completes the computation of the eigenvalues of $T$ and
therefore of the eigenvalues of $\nabla$ on $\Om_1$. We
summarize our results:

The space $\Om_1$ is the direct sum of the spaces
$\EEONE$; that is,
\[\Om_1=\bigoplus_{(M,p,k)}\EEONE\]
and the spaces $\EEONE$ are invariant under~$\nabla$. Hence,
the eigenvalue problem for $\nabla$ is reduced to the problem
of calculating the eigenvalues on these spaces of
dimension~1,2 or~3.

\smallskip
{\narrower

\noindent {\bf Case 1}: If $\dim\EEONE=1$ and $M>1$, then
$\nabla$ has eigenvalue $q^{14}\l_k^2$ or $q^6\l_k^2$ on
$\EEONE$; if $\dim\EEONE=1$ and $M=1$, then $\nabla$ has
eigenvalues ${q^{18}(1+q^2)^2}$, $q^{10}$ and ${q^2(1+q^2)^2}$
on $\E_1(1,0,-1)$, $\E_1(1,0,0)$ and $\E_1(1,0,+1)$,
respectively.

\noindent {\bf Case 2}: If $\dim\EEONE=2$, then, on $\EEONE$,
$\nabla$ has eigenvalue $\nu_k$ and either
$({q^5-q^7\l_k})^2$, if $k=-m$, or $({q^5+q^3\l_k})^2$, if
${k=+m}$.

\noindent {\bf Case 3}: If $\dim\EEONE=3$, then $\nabla$ has
eigenvalues $\nu_k$, $(\mu^+_k)^2$ and $(\mu^-_k)^2$, where
\[
\mu_k^\pm=(B_k\pm\sqrt{B_k^2-4C_k})/2,\] and
\[ B_k=(q^3-q^7)\l_k+q^5 \quad\quad{\rm and}\quad\quad
C_k= -q^{10}\nu_k.\]}
\smallskip

Because of our choice of inner product, the Hodge operator~$L$
commutes with the Laplacian~$\nabla$. Hence, the eigenvalues
of $\nabla$ on $\Om_2$ are the same as those on $\Om_1$. The
eigenvalues of $\nabla$ on $\Om_0$ are determined by the
equation $\nabla W^M_{pk}=\nu_k W^M_{pk}$. Since
$L(\Om_0)=\Om_3$, and $L$ commutes with $\nabla$, the
eigenvalues of $\nabla$ on $\Om_3$ are the same as those on
$\Om_0$.  We have therefore calculated all the eigenvalues
of~$\nabla$.

Now we turn to the question of trying to fit this
three-dimensional differential calculus into the framework of
Connes' non-commutative differential geometry. For this we
shall need to obtain some asymptotic estimates for the
eigenvalues we have just computed.

\begin{thm} \label{thm: estimates theorem} There exists a
positive constant $C(q)$, depending only on~$q$ and not on
$M$, $p$ and $k$, such that $\nabla\ge C(q)\max(q^{-8k},1)$ on
$\EEONE$.
\end{thm}

\demo{Proof} Since $\EEONE$ admits an orthonormal basis of
eigenvectors of $\nabla$, we need only show that if $\mu$ is
an eigenvalue of $\nabla$ on $\EEONE$, then $\mu\ge
C(q)\max(q^{-8k},1)$. We shall prove this only in the case
that $\dim(\EEONE)=3$. The cases where $\dim(\EEONE)$ is equal
to~2 or~1 can be dealt with similarly, but much more simply.

We first need to establish some lower bounds for $\l_k^2$,
when $k\ne 0$. If $k>0$, then $k\ge 1/2$ and therefore
$q^{-4k}\ge q^{-2}$, from which it follows immediately that
$\l_k\ge 1$. If $k<0$, then $k\le -1/2$ and therefore
$q^{-4k}\le q^2$. Hence, $\l_k\le -q^2$. In either case,
$\l_k^2\ge q^4$, provided $k\ne 0$. Moreover,
\begin{eqnarray*}
\frac{q^{-8k}}{\lambda_k^2} & = & \frac{q^{-8k} - 2 q^{-4k}
+1}{\lambda_k^2} + 2 \, \frac{q^{-4k}-1}{\lambda_k^2} +
\frac{1}{\lambda_k^2}
\\ & = & \frac{(1-q^2)^2}{q^4} \, \frac{\lambda_k^2}{\lambda_k^2} + \frac{2
(1-q^2)}{q^2} \frac{\lambda_k}{\lambda_k^2}   +
\frac{1}{\lambda_k^2} \\ & \leq & \frac{(1-q^2)^2}{q^4} +
\frac{2 (1-q^2)}{q^2}\, \frac{1}{|\lambda_k|} +
\frac{1}{\lambda_k^2} \\ & \leq & \frac{(1-q^2)^2}{q^4} +
\frac{2 (1-q^2)}{q^4} +  \frac{1}{q^4} = \frac{(2-q^2)^2}{q^4}
\ .
\end{eqnarray*}

Since $q^4\ge q^4(2-q^2)^{-2}$, we have $\l_k^2 \ge
q^4(2-q^2)^{-2}\max(q^{-8k},1)$, provided $k\ne 0$.

Since $\nu_k\ge \l_k^2$, by Equation~(\ref{eqn: nu-eigenvalue
equation}), this gives lower bounds for $\nu_k$ also, when
$k\ne 0$. If $k=0$, we use the inequality $\nu_0\ge q^2c_0^2$.
Since $\dim\EEONE=3$, we have $m\ge 1$ and therefore
$q^{2m}\le q^2$ and $q^{-2m}\ge q^{-2}$. Hence, we have
$c_0^2\ge
{q^2(1-q^{2m})(q^{-2}-q^2)(1-q^2)^{-2}}={(1-q^{2m})(1+q^2)(1-q^2)^{-1}}
\ge {1+q^2}$. Thus, $\nu_0\ge q^4$. Therefore, for all $k$,
\[\nu_k\ge a(q)\max(q^{-8k},1), \quad {\rm where} \quad
a(q)=q^4(2-q^2)^{-2}.\]

We turn now to the question of obtaining similar inequalities
for the other two eigenvalues of $\nabla$ on $\EEONE$. Recall
that these are $(\mu_k^+)^2$ and $(\mu_k^-)^2$, where
$\mu_k^\pm$ are given by Equation~(\ref{eqn: mu-eigenvalue
equation}). Hence, if $\mu$ is one of $\mu^+$ or $\mu^-$, then
\[ \mu^2=(X_k\pm\sqrt{X_k^2-4C^2_k})/2,\]
where we set $X_k=B_k^2-2C_k$. It follows from
Equations~(\ref{eqn: B-equation}) and~(\ref{eqn: C-equation})
that
\begin{equation} X_k=
q^6[(1+q^8)\l_k^2+q^4+2q^8c_{k+1}^2+2q^2c_k^2].\end{equation}
In particular, $X_k\ne 0$, since $X_k\ge q^{10}$. Note that
$\mu^2\ge
(X_k-\sqrt{X_k^2-4C^2_k})/2=2C_k^2(X_k+\sqrt{X_k^2-4C_k^2})^{-1}\ge
C_k^2/X_k$, because $\sqrt{X^2_k-4C_k^2}\le X_k$. Since
$C_k=-q^{10}\nu_k$ and $\nu_k\ne 0$, therefore $C_k\ne 0$.
Also, since $C_k^2=q^{20}\nu_k^2$, we get $C_k^2$ is greater
than or equal to each of the numbers $q^{20}\l_k^4$,
$q^{20}\l_k^2c_{k+1}^2$ and $q^{22}\l_k^2c_k^2$, by
Equation~(\ref{eqn: nu-eigenvalue equation}). Moreover,
$C^2_k\ge q^{20}a(q)\nu_k\ge q^{20}a(q)^2q^{-8k}$. Using these
inequalities and the inequality $q^{-8k}/\l_k^2\le 1/a(q)$,
where $k\ne 0$, that we established above, we get
\begin{eqnarray*}
\frac{q^{-8k}}{\mu^2}&\le& \frac{X_kq^{-8k}}{C_k^{2}}=
\frac{q^{-8k}q^6[(1+q^8)\l_k^2+q^4+2q^8c^2_{k+1}+2q^2c_k^2]}{C_k^2}
\\ &\le&
q^6[(1+q^8)\frac{q^{-8k}}{q^{20}\l_k^2}+q^4\frac{q^{-8k}}{C_k^2}+
2q^8\frac{q^{-8k}}{q^{20}\l_k^2}+2q^2\frac{q^{-8k}}{q^{22}\l_k^2}]
\\ &\le&
q^6[(1+q^8)q^{-20}a(q)^{-1}+q^{-16}a(q)^{-2}+2q^{-12}a(q)^{-1}+ 2q^{-20}a(q)^{-1}].
\end{eqnarray*}

If we denote the reciprocal of the righthand side of this last
inequality by $b_1(q)$, then we have shown that $\mu^2\ge
b_1(q)q^{-8k}$, whenever $k\ne 0$.

If we leave out the factor $q^{-8k}$ in the above inequalities
and use the fact that, for $k\ne 0$, we have $\l_k^2\ge q^4$
and therefore $C_k^2\ge q^{28}$, we get
\begin{eqnarray*} \frac{1}{\mu^2} &\le &
q^6[(1+q^8)\frac{1}{q^{20}\l_k^2}+
\frac{q^4}{C_k^2}+2q^8\frac{1}{q^{20}\l_k^2}+2q^2\frac{1}{q^{22}\l_k^2}]
\\ &\le &
q^6[(1+q^8)\frac{1}{q^{20}q^4}+q^4\frac{1}{q^{28}}+
2q^8\frac{1}{q^{20}q^4}+2q^2\frac{1}{q^{22}q^4}].\end{eqnarray*}
Hence, if $b_2(q)$ denotes the reciprocal of the righthand
side of the last inequality, we have shown that $\mu^2\ge
b_2(q)$.

Finally, we have to consider the case where $k=0$. (Then
$\l_k=0$ and we cannot use the preceding estimates.) However,
in this case, we  have
$C_0^2=q^{20}\nu_0^2={q^{20}(c_1^2+q^2c_0^2)^2}$, by
Equation~(\ref{eqn: nu-eigenvalue equation}); hence, $C_0^2\ge
q^{20}c_1^4\ge q^{20}c_1^2$ and $C_0^2\ge q^{24}c_0^4\ge
q^{26}c_0^2$. Also, since $\nu_0\ge a(q)$, we have $C_0^2\ge
q^{20}a(q)^2$. It follows from these inequalities that
\begin{eqnarray*} \frac{1}{\mu^2}\le \frac{X_0}{C_0^2}=
q^6[\frac{q^4}{C_0^2}+\frac{2q^8c_1^2}{C^2_0}+\frac{2q^2c_0^2}{C_0^2}]\le
q^6[\frac{q^4}{q^{20}a(q)^2}+\frac{2q^8}{q^{20}}+
\frac{2q^2}{q^{26}}].\end{eqnarray*} Thus, $\mu^2\ge b_3(q)$,
where $b_3(q)$ is the reciprocal of the righthand side of the
last inequality.

Set $b(q)=\min(b_1(q),b_2(q),b_3(q))$. Then $b(q)$ is a
positive constant depending only on $q$ and not on $M$, $p$
and $k$ and, for $\mu=\mu^\pm_k$ we have shown that $\mu^2\ge
b(q)\max(q^{-8k},1)$, for all $k$. To complete the proof (in
the case that $\dim(\EEONE)=3$), take $C(q)=\min(a(q),b(q))$.
\qed

We now define a decomposition for the space $\Om$ that is
analogous to that given for $\Om_1$: If $m$ is a non-negative
integer or half-integer, $M=2m+1$ and if
$p,k={-m,-m+1,\dots,+m}$, set $\E_0(M,p,k)=\bC W^M_{pk}$; if
$k=-m-1$ or $k=m+1$, set $\E_0(M,p,k)=0$. Set
$\E_2(M,p,k)=L(\EEONE)$ and $\E_3(M,p,k)=L(\E_0(M,p,k))$. We
have
\[ \Om=\bigoplus_{(M,p,k)} \EE \ \ , \]
where $\EE$ denotes the orthogonal sum $\oplus_{i=0}^3 \,\E_i(M,p,k)$. It follows from Theorem~\ref{thm: estimates
theorem} that there exists a positive constant $C(q)$, depending only on~$q$ and not on $M$, $p$ or $k$, such that, for
$M>1$, the restriction $\nabla_\E$ of $\nabla$ to $\EE$ satisfies the following inequalities
\begin{equation}\label{eqn: 3-estimates}\nabla_{\E}\ge C(q),\quad
\nabla_{\E}\ge C(q)q^{-8k}, \quad \nabla_{\E}\ge
C(q)q^{-4k}.\end{equation} For $M=1$, these inequalities hold
for $\nabla_\E$ equal to the restriction of $\nabla$ to
${\E=\E_1(1,p,k)\oplus\E_2(1,p,k)}$ ($\nabla=0$ on
$\E_0(1,p,k)\oplus\E_3(1,p,k)$). The third of the inequalities
in~(\ref{eqn: 3-estimates}) is an easy consequence of the
first two. The required inequalities for $\E_0(M,p,k)$, when
this space is non-zero, and $M>1$, follow from the fact that
the eigenvalue of $\nabla$ on $\E_0(M,p,k)=\bC W^M_{pk}$ is
$\nu_k$, which is also an eigenvalue of $\nabla$ on the
corresponding space $\EEONE$. We are also using the fact that
the eigenvalues of $\nabla$ on $\E_2(M,p,k)$ and $\E_3(M,p,k)$
are the same as those on $\E_1(M,p,k)$ and $\E_0(M,p,k)$,
respectively.

We set \[\GG(M,k)=\oplus_{p=-m}^{p=+m}\EE \quad{\rm for}\quad k={-m-1,-m,\dots,m,m+1}.\] Then $\Om=\oplus_{(M,k)}
\GG(M,k)$, as an orthogonal sum by the orthogonality properties of the $W^M$.

Since $\GG(M,k)$ is finite-dimensional, and $\nabla$ is
self-adjoint, the restriction $\nabla_{(M,k)}$ to $\GG(M,k)$
is diagonalisable;  the inequalities in~(\ref{eqn:
3-estimates}), for $M>1$, imply corresponding inequalities for
the eigenvalues of $\nabla_{(M,k)}$, from which it follows
that \begin{equation}\label{eqn: alternative 3-estimates}
\nabla_{(M,k)}\ge C(q),\quad \nabla_{(M,k)}\ge C(q)q^{-8k},
\quad \nabla_{(M,k)}\ge C(q)q^{-4k}.\end{equation}

Now define a bounded operator $R$ on $\Om$ by setting
$R=(1+\nabla)^{-1/2}$. This is to be understood more precisely
as saying that $R$ is the direct sum of the operators
$R_{(M,k)}$, where $R_{(M,k)}=(1+\nabla_{(M,k)})^{-1/2}$.
Since $R_{(M,k)}\le 1$, we have $\n R\le 1$. However, we can
do better: it follows easily from the inequalities
in~(\ref{eqn: alternative 3-estimates}) that
\begin{equation} \label{eqn: norm estimates} \n{R_{(M,k)}}\le
C(q)^{-1/2}\min(q^{2k},q^{4k}).\end{equation} (We can easily
modify $C(q)$ by making it smaller so that (\ref{eqn: norm
estimates}) also holds for $M=1$.)

We recall from \cite[Equation~2.16]{WOR1} that there are
linear functionals $f_1,f_2,f_3$ on $\sa_q$ such that $\eta_r
a=(\id \otimes f_r)\D(a) \eta_r$, for all $a\in \sa_q$.
Moreover, $f_3=f_1$ and $(f_1\otimes
f_1)\D=f_2=\e+(q^{-2}-1)\chi_2$, where $\e$ is the co-unit of
$\sa_q$.

Set $B_r=f_r(W^M)$. Then
$(B_2)_{pk}=\delta_{pk}+(q^{-2}-1)\l_k\delta_{pk}=q^{-4k}\delta_{pk}$.
Also, $B_1=B_3$ and $B_1^2=B_2$. Since $B_2$ is diagonal, with
distinct diagonal entries, and $B_1$ commutes with $B_2$,
therefore $B_1$ is also diagonal. In fact, the diagonal
entries of $B_1$ are positive~\cite[Equation~5.53]{WOR1}, so
$(B_1)_{pk}=q^{-2k}\d_{pk}$. Hence, for all~$r$,
\begin{equation}\label{eqn: commutation equations} \eta_r
W^M_{pk}=(\id\otimes
f_r)\D(W^M_{pk})\eta_r=(W^MB_r)_{pk}\eta_r
=q^{-2\e_rk}W^M_{pk}\eta_r,\end{equation} where $\e_1=\e_3=1$
and $\e_2=2$.

Suppose that $M$ is an integer, $M=2m+1$, and that $k$ is a
half-integer or integer. We define $\GG(M,k)=0$, if $M\le 0$
or if $k$ does not belong to the set
${\{-m-1,-m,\dots,m,m+1\}}$.

\begin{lem} \label{lem: G-lemma} There exists a  positive number $C(q)$
(depending on $q$ but not on $M$ or $k$) such that, for $M\ge
1$ and ${k=-m-1,-m,\dots,m,m+1}$, we have
\begin{enumerate}
\item $M_1(\GG(M,k)) \subseteq \GG(M,k-1)$ and $\n{(M_1)_{\GG(M,k)}} \le
C(q)q^{-2k}$;
\item $M_2(\GG(M,k)) \subseteq \GG(M,k)$ and $\n{(M_2)_{\GG(M,k)}}\le
C(q)q^{-4k}$;
\item $M_3(\GG(M,k)) \subseteq \GG(M,k+1)$ and $\n{(M_3)_{\GG(M,k)}}\le
C(q)q^{-2k}$.
\end{enumerate}
\end{lem}

\demo{Proof} The subset inclusions $M_1(\EE) \subseteq
\E(M,p,k-1)$, $M_2(\EE) \subseteq \E(M,p,k)$ and $M_3(\EE)
\subseteq \E(M,p,k+1)$ are easy consequences of Equation
Set~(\ref{eqn: commutation equations}). The subset inclusions
in Conditions~(1)--(3) follow immediately. The space $\EE$ is
the linear span of the orthogonal vectors
\begin{eqnarray*}
\begin{array}{cccc} &g^p_1& \\ g^p_2&g^p_3&g^p_4\\
g^p_5&g^p_6&g^p_7\\ &g^p_8
&\end{array}=\begin{array}{cccc}&W^M_{pk}&\\
\ W^M_{p,k+1}\eta_1& W^M_{pk}\eta_2&W^M_{p,k-1}\eta_3\\
W^M_{p,k+1}L(\eta_1)& W^M_{pk}L(\eta_2)&W^M_{p,k-1}L(\eta_3)\\
&W^M_{pk}L(1)&\end{array}\\ =
\begin{array}{cccc}&W^M_{pk}&\\
W^M_{p,k+1}\eta_1& W^M_{pk}\eta_2&W^M_{p,k-1}\eta_3\\
q^{-7}W^M_{p,k+1}\eta_1\eta_2&
-q^{-6}W^M_{pk}\eta_3\eta_1&qW^M_{p,k-1}\eta_2\eta_3\\
&W^M_{pk}\ta.&\end{array}
\end{eqnarray*}
Hence, setting $b_k=(B_1)_{kk}=q^{-2k}$ and $b_{-m-1}=0$, and
using Equation Set~(\ref{eqn: commutation equations}) and the
commutation relations for the products of the elements
$\eta_1$, $\eta_2$ and $\eta_3$, we get that the images of the
above vectors under $M_1$ are given by
\begin{eqnarray*}
\begin{array}{cccc}&b_kW^M_{pk}\eta_1&\\
0&
-q^{-4}b_kW^M_{pk}\eta_2\eta_1&-q^{-2}b_{k-1}W^M_{p,k-1}\eta_3\eta_1\\
0& 0 &q^{-5}b_{k-1}W^M_{p,k-1}\eta_2\eta_3\eta_1\\
&0&\end{array}
\end{eqnarray*}
that is, by
\begin{eqnarray*}
\begin{array}{cccc} &b_kg^p_1\eta_1& \\0&-q^{-4}b_kg^p_3\eta_1
&-q^{-2}b_{k-1}g^p_4\eta_1\\ 0&0&q^{-6}b_{k-1}g^p_7\eta_1\\ &0
&\end{array}.\end{eqnarray*}

Now let $Z_r$ be the linear span of the vectors
${g^{-m}_r,g^{-m+1}_r,\dots,g^m_r}$. Note that the linear map,
${R_{\eta_1}\colon\om\mapsto \om\eta_1}$, is bounded, since it
is the tensor product of the identity map on $\sa_q$ and its
own restriction to $\ominv$, and the latter is bounded since
$\ominv$ is finite-dimensional. The norm of $R_{\eta_1}$ clear
depends only on $q$ (provided we fix, as we have, the inner
product on $\ominv$). Then, if ${\om\in Z_1}$,
$M_1(\om)=b_k\om\eta_1$, so that $\n{(M_1)_{Z_1}}\le
b_k\n{R_{\eta_1}}$. Similarly, $\n{(M_1)_{Z_3}}\le
q^{-4}b_k\n{R_{\eta_1}}$, $\n{(M_1)_{Z_4}}\le
q^{-2}b_{k-1}\n{R_{\eta_1}}$, $\n{(M_1)_{Z_7}}\le
q^{-6}b_{k-1}\n{R_{\eta_1}}$ and $M_1=0$ on $Z_2$, $Z_5$,
$Z_6$ and $Z_8$. Now $\GG(M,k)$ is the orthogonal direct sum
of the spaces $Z_r$, that is, $\GG(M,k)=\oplus_{r=1}^8 Z_r $.
It follows that
\begin{eqnarray*}
\n{(M_1)_{\GG(M,k)}}&\!=\!&\max_r(\n{(M_1)_{Z_r}})\le
\max(b_k,q^{-4}b_k,q^{-2}b_{k-1},q^{-6}b_{k-1})\n{R_{\eta_1}}\\
&\!\le\!&\max(q^{-2k},q^{-4}q^{-2k},q^{-2k},q^{-4}q^{-2k})\n{R_{\eta_1}}.
\end{eqnarray*}
Setting $a_1(q)=\max(1,q^{-4})\n{R_{\eta_1}}$, we get
$\n{(M_1)_{\GG(M,k)}}\le a_1(q)q^{-2k}$.

One can obtain, by similar methods, positive numbers $a_2(q)$
and $a_3(q)$ such that $\n{(M_2)_{\GG(M,k)}}\le a_2(q)q^{-4k}$
and $\n{(M_3)_{\GG(M,k)}}\le a_3(q)q^{-2k}$. Setting
$C(q)=\max(a_1(q),a_2(q),a_3(q))$, we get the inequalities in
Conditions~(1)--(3). \qed

If $a\in \sa_q$, we denote by $L_a$ the bounded operator on
$\Om$ obtained by left multiplication by~$a$.

\begin{thm} \label{thm:Jaffe} Let $\b$ and $\d$ be nonnegative numbers for which
${\b+\d=1}$. Then $R^\b[L_a,D]R^\d$ is a bounded operator on
$\Om$, for all $a\in \sa_q$. \end{thm}

\demo{Proof} We shall prove that $R^\b[L_a,d]R^\d$ is bounded,
for all $a\in \sa_q$. This will suffice, since then, replacing
$a$ by $a^*$ and interchanging the roles of $\b$ and $\d$, and
using the fact that $L_a^*=L_{a^*}$, we get the adjoint
$-R^\b[L_a,d^*]R^\d$ of $R^\d[L_{a^*},d]R^\b$ is bounded.
Hence, $R^\b[L_a,D]R^\d$ is bounded.

Since $[L_a,d](\om)=-(d a)\om$, for all $\om\in \Om$, we have
$[L_a,d]=-\sum_{r=1}^3L_{\chi_r*a}M_r$. Consequently, it
suffices to show that $R^\b L_aM_rR^\d$ is bounded, for all
$a\in \sa_q$ and $r=1,2,3$.

Let $\T$ denote the set of all bounded operators in $B(\Om)$
for which there exist $N\in \Z$ and $l\in \frac 12\Z$ such
that $T\GG(M,k)\subseteq \GG(M+N,k+l)$, for all $M\in \Z$ and
${k\in \frac 12\Z}$. One easily checks that $\T$ is
self-adjoint and closed under multiplication, so that its
linear span $\T'$ is a self-adjoint subalgebra of $B(\Om)$.
Hence, by the following lemma, $L_a\in \T'$, for all $a\in
\sa_q$. Thus, to prove the theorem, we need only show now that
$R^\b TM_rR^\d$ is bounded, for all $T\in \T$. We shall show
this only in the case that $r=1$. The cases where $r=2$ or
$r=3$ have similar proofs.

Before proceeding let us note first that by combining
Inequality~(\ref{eqn: norm estimates}) and the inequality in
Condition~1 of Lemma~\ref{lem: G-lemma} we get a positive
constant $C$ depending only on $q$ such that $\n{R_{(M,k)}}\le
Cq^{2k}$ and $\n{(M_1)_{\GG(M,k)}}\le Cq^{-2k}$.

Suppose then $T\in\T$ and $N\in\Z$ and $l\in \frac 12\Z$ are
such that for all $M=1,2,\dots$ and ${k=-m-1,m,\dots,m,m+1}$,
we have $T\GG(M,k)\subseteq \GG(M+N,k+l)$. Let $\om\in
\GG(M,k)$. Then $TM_1R^\d\om\in\GG(M+N,k+l-1)$. Hence,
\begin{eqnarray*}\n{R^\b TM_1R^\d\om}^2 &\le&
C^{2\b}q^{4(k+l-1)\b}\n{TM_1R^\d\om}^2 \\ &\le &
C^{2\b}q^{4(k+l-1)\b}\n{T}^2C^2q^{-4k}\n{R^\d\om}^2 \\ &\le &
C^{2\b}q^{4(k+l-1)\b}\n{T}^2C^2q^{-4k}C^{2\d}q^{4k\d}\n{\om}^2
\\ &=& C^{2(\b+\d)}C^2\n T^2q^{4k(\b+\d)}q^{4(l-1)\b}
q^{-4k}\n\om^2 \\ &=& C^4\n
T^2q^{4(l-1)\b}\n\om^2.\end{eqnarray*} It follows now, from
orthogonality of the image spaces $R^\b TM_1R^\d\GG(M,k)$,
that $R^\b TM_1R^\d$ is bounded.~\qed

\begin{lem}  \label{lemma:927522} If $\a$ and $\g$ are the canonical generators of
$\sa_q$ (see \cite{WOR1}), then $L_\a$ and $L_\g$ belong to
the linear span of the set $\T$ defined in the preceding
proof.
\end{lem}

\medskip

Let us introduce some extra notation. We will use the elements $x_k$ introduced in the proof of
\cite[Theorem~5.8]{WOR1}, but add an extra parameter into the notation. Let $m \in \frac{1}{2} \Z$. Then we define
$x_k^m = \a^{m+k} (\gamma^*)^{m-k}$ for $k=-m,\ldots,m$.

For $k \in \Z$ and $y \in \sa_q$ we set $y(k) = y^k$ if $k \geq 0$ and $y(k) = (y^*)^{-k}$ if $k < 0$. We also define
$a(k,l) = \a(k)\,\gamma(l)$ for all $k,l \in \Z$. It follows from~ \cite[Theorem~1.2]{WOR1} that the family
$\{\,a(k,l)\,(\gamma^* \gamma)^m \mid k,l \in \Z, m \in \N \, \}$ forms a linear basis for $\sa_q$. So if
$\sa_q(\gamma^* \gamma)$ denotes the unital $^*$-subalgebra of $\sa_q$ generated by $\gamma^* \gamma$, we see that
$\sa_q$ is the direct sum $\oplus_{u,v \in \Z}  a(u,v) \sa_q(\gamma^* \gamma)$.

\begin{lem} \label{lem:7219201}
Consider $M \in \N_0$ and $m \in \frac{1}{2}\,\Z$ so that $M = 2m+1$. If $j,k \in \{-m,\ldots,m\}$, then  $W^M_{jk}
\in a(j+k,k-j)\,\sa_q(\gamma^* \gamma)$
\end{lem}
\begin{proof} We proceed by induction on $M$. Since $W^1_{00} = 1$, the lemma is certainly true for $M = 1$. Next we
suppose that $M > 1$ and that the lemma is true for $M-1$.
\begin{trivlist}
\item[\ \,(1)] Let $j \in \{-m,\ldots,m\}$. Clearly, $x_{-m}^m = (\gamma^*)^{2m}$ thus
$\sum_{i=-m}^m x_i^m \ot W^M_{i,-m} = \D(\gamma^*)^{2m} = (\a \ot \gamma^*  + \gamma^* \ot \a^*)^{2m}$. Since $\a \ot
\gamma^*$ and $\gamma^* \ot \a^*$ commute up to a scalar, there exist complex numbers $r_{-m},\ldots,r_m$ such that
\begin{eqnarray*}
(\gamma^* \ot \a^* + \a \ot \gamma^*)^{2m} & = & \sum_{i=-m}^m r_i\,(\a \ot \gamma^*)^{i+m} (\gamma^* \ot \a^*)^{m-i}
\\ & = & \sum_{i=-m}^m r_i\,q^{(m+i)(m-i)}\,x_i^m \ot a(i-m,-i-m)\ ;
\end{eqnarray*}
thus, $W^M_{j,-m} = r_j \, q^{(m+j)(m-j)}\,a(j-m,-j-m)$.
\medskip
\item[\ \,(2)] Let $k \in \{-m+1,\ldots,m\}$ and $j \in \{-m,\ldots,m\}$.  Clearly,
$\a\,x^{m-\frac{1}{2}}_{k-\frac{1}{2}} = x^m_k$. Hence
\begin{eqnarray*}
& & \sum_{i=-m}^m x^m_i \ot W^M_{ik} = \D(x^m_k) = \D(\a\,x^{m-\frac{1}{2}}_{k-\frac{1}{2}}) =
\D(\a)\,\D(x^{m-\frac{1}{2}}_{k-\frac{1}{2}})
\\ & & \hspace{3ex} = (\a \ot \a -q\,\gamma^* \ot
\gamma)\,\,\sum_{i=-m+\frac{1}{2}}^{m-\frac{1}{2}}\,x^{m-\frac{1}{2}}_i \ot W^{M-1}_{i,k-\frac{1}{2}}
\\ & & \hspace{3ex} = \sum_{i=-m+\frac{1}{2}}^{m-\frac{1}{2}}\,x^{m}_{i+\frac{1}{2}} \ot \a\,W^{M-1}_{i,k-\frac{1}{2}}
- \sum_{i=-m+\frac{1}{2}}^{m-\frac{1}{2}}\,q^{-m-i+\frac{3}{2}}\,     x^m_{i-\frac{1}{2}} \ot
\gamma\,W^{M-1}_{i,k-\frac{1}{2}}
\\ & & \hspace{3ex} = \sum_{i=-m+1}^m\,x^{m}_i \ot \a\,W^{M-1}_{i-\frac{1}{2},k-\frac{1}{2}}
- \sum_{i=-m}^{m-1}\,q^{-m-i+1}\,     x^m_i \ot \gamma\,W^{M-1}_{i+\frac{1}{2},k-\frac{1}{2}}
\end{eqnarray*}
This implies the existence of complex numbers $c,d$ such that $W^M_{jk}$ equals \newline $c\,
\a\,W^{M-1}_{j-\frac{1}{2},k-\frac{1}{2}} + d \, \gamma\,W^{M-1}_{j+\frac{1}{2},k-\frac{1}{2}}$. Therefore the
induction hypothesis implies that
$$W^M_{jk} \in \a\,a(j+k-1,k-j)\, \sa_q(\gamma^* \gamma) +  \gamma\,a(j+k,k-j-1)\,\sa_q(\gamma^* \gamma)
$$  Since $\a\,a(p,q) \in a(p+1,q)\,\sa_q(\gamma^* \gamma)$ and
$\gamma\,a(p,q) \in a(p,q+1)\,\sa_q(\gamma^* \gamma)$, it follows that $W^M_{jk} \in a(j+k,k-j) \, \sa_q(\gamma^*
\gamma)$.
\end{trivlist}
So we have proven that $W^M_{jk} \in a(j+k,k-j)\,\sa_q(\gamma^* \gamma)$ for all possible values of $j$ and $k$.
\end{proof}

Lemma \ref{lemma:927522} will be an immediate consequence of the next one.

\begin{lem}
Consider $M \in \N_0$ and $m \in \frac{1}{2}\,\Z$ so that $M = 2m+1$. If $j,k \in \{-m,\ldots,m\}$, then $\a\,W^M_{jk}
\in  \C\,W^{M-1}_{j+\frac{1}{2},k+\frac{1}{2}} + \C\,W^{M+1}_{j+\frac{1}{2},k+\frac{1}{2}}$ and $\gamma\,W^M_{jk} \in
\C\,W^{M-1}_{j-\frac{1}{2},k+\frac{1}{2}} + \C\,W^{M+1}_{j-\frac{1}{2},k+\frac{1}{2}}$.
\end{lem}
\begin{proof}
One easily checks that $W^2 = \left(\begin{array}{cc} \a^* & -q \, \gamma \\ \gamma^* & \a \end{array}\right)$. Since
$W^1_{00} = 1$, this implies that the lemma is certainly true if $M = 1$. From now on, we suppose that $M > 1$.

By \cite[Theorem~5.11]{WOR1}, we know that the tensor product representation $W^2 \ot W^M$ is equivalent to $W^{M-1}
\oplus W^{M+1}$, that is, there exists an invertible matrix $Q \in M_{2M}(\C)$ such that $W^2 \ot W^M = Q(W^{M-1}
\oplus W^{M+1}) Q^{-1}$. But $\a \, W^M_{jk}$ appears as a matrix element of $W^2 \ot W^M$. It follows that
$\a\,W^M_{jk}$ belongs to $\langle \, W^{M-1}_{rs} \mid r,s=-m+\frac{1}{2},\ldots,m-\frac{1}{2} \,\rangle + \langle \,
W^{M+1}_{rs} \mid r,s=-m-\frac{1}{2},\ldots,m+\frac{1}{2} \, \rangle$.

By the previous lemma,  $\a \, W^M_{jk} \in a(j+k+1,k-j)\,\sa_q(\gamma^* \gamma)$. On the other hand, we know that
$W^{M-1}_{rs} \in a(r+s,s-r)\,\sa_q(\gamma^* \gamma)$ and $W^{M+1}_{r's'} \in a(r'+s',s'-r')\,\sa_q(\gamma^* \gamma)$,
for all $r,s \in \{-m+\frac{1}{2},\ldots,m-\frac{1}{2}\}$ and  $r',s' \in \{-m-\frac{1}{2},\ldots,m+\frac{1}{2}\}$.
Note also that $j+k +1 = r+s$ and $k-j = s-r$ $\Leftrightarrow$ $r=j+\frac{1}{2}$ and $s=k+\frac{1}{2}$. Since $\sa_q$
is the direct sum $\oplus_{u,v} a(u,v)\,\sa_q(\gamma^* \gamma)$, this implies that $\a\,W^M_{jk} \in
\C\,W^{M-1}_{j+\frac{1}{2},k+\frac{1}{2}} + \C\,W^{M+1}_{j+\frac{1}{2},k+\frac{1}{2}}$. The statement concerning
$\gamma$ is proven in a similar way.
\end{proof}

We used the next result earlier in this section but prove it here because it is an easy consequence of Lemma
\ref{lem:7219201}.

\begin{lem} \label{lem:orthogonality}
Consider $M \in \N_0$ and $m \in \frac{1}{2}\,\Z$ so that $M = 2m+1$. If $p,p',j,k \in \{-m,\ldots,m\}$ and $(p,j)
\not= (p',k)$, then $h((W^M_{pj})^* W^M_{p'k}) = 0$.
\end{lem}
\begin{proof}
If $i,j,i',j' \in \Z$ and $b,c \in \sa_q(\gamma^* \gamma)$, equation \cite[A~1.8]{WOR2} implies that $h((a(i',j')\,c)^*
(a(i,j)\,b)) = 0$ if  $(i,j) \not= (i',j')$. By assumption, $(p+j,j-p) \not= (p'+k,k-p')$. As a consequence, lemma
\ref{lem:7219201} implies that $h((W^M_{pj})^* W^M_{p'k}) = 0$.
\end{proof}

\medskip

It is very unlikely that the Dirac operator considered here, fits into the framework of Connes' non-commutative
geometry (although we do not have a proof of this fact).  This provides another indication that generalisations of this
non-commutative geometry have to be studied. One such generalisation is considered in \cite{Jaffe}, but we could not
prove that our Dirac operator fits into this more general framework for non-commutative geometry. The problem lies in
the fact that we can only prove Theorem \ref{thm:Jaffe} under the assumption that $\beta + \delta = 1$ whereas the
theory in \cite{Jaffe} relies on the fact that this same theorem is true if $\beta + \delta < 1$ (see \cite[Section
V.5]{Jaffe}).


In the next part we collect some concrete formulas for the one-parameter groups and twist automorphism. Let $\rho$ the
modular group and $\tau$ denote the scaling group of $\sa_q$. It follows from \cite[Appendix~A1]{WOR2} that $\rho_z(\a)
= q^{-2iz}\,\a$, $\rho_z(\gamma) = \gamma$, $\tau_z(\a) = \a$ and $\tau_z(\gamma) = q^{2iz}\, \gamma$, for all $z \in
\C$.

We denote the twist automorphism of $\int$ by $\s$. We know from \cite[Table~1]{WOR1} that $\ta\a = q^{-4}\,\a$,
$\ta\gamma = q^{-4}\,\gamma$. Thus $\s(\a) = q^{-4}\,\rho_i(\a) = q^{-2}\,\a$. Similarly, $\s(\a^*) = q^2\,\a^*$,
$\s(\gamma) = q^{-4}\, \gamma$ and $\s(\gamma^*) = q^4 \,\gamma^*$. Since $\s$ commutes with the differential $d$, it
follows from \cite[Table~3]{WOR1} that $\s(\eta_1) = q^6\,\eta_1$, $\s(\eta_2) = \eta_2$ and $\s(\eta_3) =
q^{-6}\,\eta_3$. In section \ref{sec: one-parameter reps} it turned out that it is beneficial to require the inner
product on $\ominv$ to be chosen in such a way that $\s$ is a positive operator. Since $\eta_1$, $\eta_2$ and $\eta_3$
are eigenvectors of $\s$ with different eigenvalues, this requirement implies that $\eta_1$, $\eta_2$, $\eta_3$ are
necessarily orthogonal with respect to such an inner product (which is true for the inner product used in this
section). Similar remarks apply to the second order forms.

By the universality of $\sa_q$ there exists a unique algebra homomorphism $\phi : \sa_q \rightarrow \sa_q$ such that
$\phi(\a) = q^{-1}\,\a$, $\phi(\a^*) = q \a^*$, $\phi(\gamma) = q^{-2}\,\gamma$ and $\phi(\gamma^*) = q^2\,\gamma^*$.
One easily checks that $\phi^2(a) = \s(a)$, $\phi(\phi(a^*)^*) = a$  and $(\tau_{\frac{i}{2}} \ot \phi) \D(a) =
\D(\phi(a))$, for all $a \in \sa_q$ (it is enough to check these equalities on the generators of $\sa_q$). Thus it
follows from Corollary \ref{one.cor1} that $h(a^* \s(a)) \geq 0$. By Theorem \ref{one.prop2} there exists a unique
exponential one-parameter group $\b$ on $\Om$ such that $\b_i = \s$.

One easily checks that $\b_z(\a) = q^{2iz}\,\a$ and $\b_z(\gamma) = q^{4iz}\,\gamma$, for all $z \in \C$. As before,
\cite[Table~2]{WOR1} implies that $\b_z(\eta_1) = q^{-6iz}\,\eta_1$, $\b_z(\eta_2) = \eta_2$ and $\b_z(\eta_3) =
q^{6iz}\,\eta_3$, for all $z \in \C$. Also note that the inner product used in this section satisfies the requirements
of Theorem \ref{thm: inner product theorem 2}.

\bigskip

\section{Appendix: one-parameter groups}

In this appendix we present a self-contained brief account of the material on exponential one-parameter groups of
automorphisms needed for this paper. We have included this since the material is not standard and we know of no
suitable reference for it in the literature.

Let $\b$ be a map from $\bC$ into the set of algebra
automorphisms on a $*$-algebra~$B$. Suppose that $\b_z(b)^* =
\b_{\bar{z}}(b^*)$ and $\b_{y+z} = \b_y \, \b_z$, for all $b
\in B$ and $y,z \in \C$. Then $\b$ is called a  {\em
one-parameter group} on $B$.

Note that $\b_0 = \id_B$ and  $(\b_z)^{-1} = \b_{-z}$, for
$z\in \C$.

\smallskip

One-parameter groups become even more interesting if extra
analyticity conditions are imposed. In the case of compact
quantum groups an exponential growth condition can also be
imposed.

We say that an analytic function $f : \C \rightarrow \C$ is of
{\em exponential type} if there exist positive numbers $M$ and
$r$ such that $|f(z)| \leq M \, e^{r \, |\text{Im\,}z|}$, for
all $z \in \C$.

It is easy to check that the functions of exponential type
form an algebra. It is equally clear that for any positive
number $\lambda$, the function $\C \rightarrow \C : z \mapsto
\lambda^{iz}$ is of exponential type. The translate of a
function of exponential type is another such function.

Consider a vector space  $V$ and a function $f : \C
\rightarrow V$. We say that $f$ is of {\em finite exponential
type} if there exist elements $v_1,\ldots,v_n \in V$ and
functions $f_1,\ldots,f_n$ of exponential type such that $f(z)
= \sum_{j=1}^n f_j(z) \, v_j$, for all $z \in \C$.

Let $f,g : \C \rightarrow V$ be two functions of finite
exponential type and $c$ a non-zero number on the imaginary
axis. Then, by \cite[Lemma~5.5]{WOR2}, $f = g$ if and only if
$f(n c ) = g(n c )$, for all $n \in \N_0$.

Consider a one-parameter group $\b$ on a $^*$-algebra $B$ such
that for every $b \in B$, the function $\C \rightarrow B : z
\mapsto \b_z(b)$ is of finite exponential type. Then we call
$\b$ an {\em exponential} one-parameter group on $B$.

Let $c$ be a non-zero number on the imaginary axis. Then an
exponential one-parameter group on $B$ is completely
determined by its value at $c$: If $\b$ and~$\g$ are
exponential one-parameter groups on $B$, then $\b = \gamma$ if
and only if $\b_c = \gamma_c$.

\smallskip

Another basic fact is the following: If $\t$ is a linear
functional on $B$ such that $\t \, \b_c = \t$, then $\t \,
\b_z = \t$, for all $z \in \C$. This is true because, for
every $b \in  B$, the function $\C \rightarrow \C : z \mapsto
\t(\b_z(b))$ is of exponential type and $\t(\b_{n c}(b)) =
\t(b)$, for all $n \in \N_0$.

\medskip\medskip

In order to check that an algebra automorphism on a
$^*$-algebra is induced by an exponential one-parameter group,
it is enough to check it on a set of generators:

\begin{thm} \label{one.prop1}
Let $B$ be a $^*$-algebra and $S$ a subset such that $S \cup
\{1\}$ generates~$B$ as a $*$-algebra. Let $\phi$ be an
algebra homomorphism on $B$ satisfying ${\phi(\phi(b^*)^*) =
b}$, for all $b \in B$. Assume for every $b \in S$ the
existence of a function $f : \C \rightarrow B$ of finite
exponential type such that $f(n i) = \phi^n(b)$, for all $n
\in \N$. Then there exists a unique exponential one-parameter
group $\b$ on $B$ such that $\b_i = \phi$.
\end{thm}

\demo{Proof} Define $T$ to be the set of elements $b \in B$
such that there exist a function $f : \C \rightarrow B$ of
finite exponential type satisfying $f(n i) = \phi^n(b)$, for
all $n \in \N$. Note that  $f$ is unique; we denote it by
$f_b$. Clearly, $1 \in T$ and $f_1(z) = 1$, for all $z \in
\C$. By assumption $S \subseteq T$.

The set $T$ is self-adjoint: To see this, let $b \in T$ and $m
\in \N_0$. Then the functions $\C \rightarrow B : z \mapsto
\phi^m(f_b(z))$ and $\C \rightarrow B : z \mapsto f_b(z+mi)$
are
 of finite exponential type and agree on the set $\N_0 \, i$.
So they must be equal on the whole complex plane. Hence, $\phi^m(f_b(-m i)) = f_b(0) = b$. Therefore, $f_b(-m i) =
\phi^{-m}(b)$ and $f_b(-m i)^* = \phi^m(b^*)$. It follows that $b^* \in T$ and $f_{b^*}(z) = f_b(\bar{z})^*$, for all
$z \in \C$.

If $b,b' \in T$ and $c,d \in \C$, linearity of each $\phi^n$
implies easily that $c b + d b' \in T$ and $f_{c b + d b'} = c
\, f_b + d \, f_{b'}$.  Also, multiplicativity of each
$\phi^n$ implies that $b b' \in T$ and $f_{b b'} = f_b \,
f_{b'}$. It follows that $T$ is a $^*$-subalgebra of $B$
containing $S$ and $1$, and therefore it is equal to $B$.

Define the function $\b$ from $\C$ into the set of mappings on
$B$ by setting $\b_z(b) = f_b(z)$, for all $z \in \C$ and $b
\in B$. The preceding considerations imply, for every $z \in
\C$, that $\b_z$ is an algebra homomorphism and that
$\b_z(b)^* = \b_{\bar{z}}(b^*)$, for all $b \in B$. It is also
clear that $\b_{ni} = \phi^n$, for all $n \in \N_0$.

Let $b \in B$, $y \in \C$, and $m,n \in \N_0$. The functions
$\C \rightarrow B : z \mapsto \b_z(\b_{mi}(b))$ and $\C
\rightarrow B : z \mapsto \b_{z+mi}(b)$ are of finite
exponential type and $\b_{ni}(\b_{mi}(b)) = \phi^n(\phi^m(b))
= \phi^{n+m}(b) = \b_{n i + m i}(b)$. Hence, both functions
agree on the whole complex plane; in particular
$\b_y(\b_{mi}(b)) = \b_{y+mi}(b)$.

Since  the functions $\C \rightarrow B : z \mapsto \b_y(\b_z(b))$ and $\C \rightarrow B : z \mapsto \b_{y+z}(b)$ are
both of finite exponential type, the preceding considerations imply that both functions agree on the whole complex
plane; hence, $\b_y(\b_z(b)) = \b_{y+z}(b)$, for all $y,z \in \C$. From all this we conclude that $\b$ is an
exponential one-parameter group on $B$. \qed

\begin{lem} \label{one.lem1}
Let $B$ be a $^*$-algebra, $\pi : B \rightarrow B$ a linear
map and $\b$ an exponential one-parameter group on $B$. Let
$c$ be a non-zero number on the imaginary axis. If $\pi\,\b_c
= \b_c \, \pi$ (respectively, $\pi\b_c=\b_{-c}\pi$), then $\pi
\, \b_z = \b_z \, \pi$ (respectively, $\pi\b_z=\b_{-z}\pi$),
for all $z \in \C$.
\end{lem}
\demo{Proof} We prove the result only in the first case, where
$\pi\,\b_c = \b_c \, \pi$; the proof of the other case is
almost the same. Let $b \in B$. It is clear that the functions
$\C \rightarrow B : z \mapsto \pi(\b_z(b))$ and $\C
\rightarrow B : z \mapsto \b_z(\pi(b))$ are of finite
exponential type. Our assumption implies that $\pi \, \b_{n c}
= \b_{nc} \, \pi$, for all $n \in \N_0$, implying that the
above functions are equal on $\N_0 \, i$, and therefore on the
whole complex plane.~\qed

\vspace{2cm}

{\noindent \bf Adresses}

\medskip

\begin{trivlist}
\item Johan Kustermans,
ITS - Basiseenheid Analyse, TU Delft, The Netherlands,
\vspace{-1ex}\begin{verbatim}j.kustermans@its.tudelft.nl\end{verbatim}
\smallskip\smallskip
\item Gerard Murphy,
Department of Mathematics, National University of Ireland, Cork, Ireland,
\vspace{-1ex}\begin{verbatim}gjm@ucc.ie\end{verbatim}
\smallskip\smallskip
\item Lars Tuset,
Faculty of Engineering, University College, Oslo, Norway,
\vspace{-1ex}\begin{verbatim}Lars.Tuset@iu.hio.no\end{verbatim}
\end{trivlist}
\end{document}